\documentclass[a4paper]{article}

\usepackage{bbm}
\usepackage{comment}

\usepackage[all]{xy}\usepackage[latin1]{inputenc}        
\usepackage[dvips]{graphics,graphicx}
\usepackage{amsfonts,amssymb,amsmath,xcolor,mathrsfs, amstext}
\usepackage{amsbsy, amsopn, amscd, amsxtra, amsthm,authblk, enumerate, dsfont}
\usepackage{upref}
\usepackage{geometry}
\usepackage[displaymath]{lineno}
\usepackage{float}

\usepackage[colorlinks,
            linkcolor=blue,
            anchorcolor=green,
            citecolor=blue
            ]{hyperref}

\numberwithin{equation}{section}

\def\e{\varepsilon}
\def\epsilon{\varepsilon}
\def\eps{\varepsilon}

\newcommand{\ol}{\overline}
\newcommand{\wt}{\widetilde}

\def\alb#1\ale{\begin{align*}#1\end{align*}}
\newcommand{\eqb}{\begin{equation}}
\newcommand{\eqe}{\end{equation}}

\DeclareMathOperator{\var}{var}

\newcommand{\bbC}{\mathbb{C}}
\newcommand{\bbD}{\mathbb{D}}
\newcommand{\bbE}{\mathbb{E}}
\newcommand{\bbH}{\mathbb{H}}

\newcommand{\bbR}{\mathbb{R}}
\newcommand{\bbP}{\mathbb{P}}
\newcommand{\bbZ}{\mathbb{Z}}

\newcommand{\cC}{\mathcal{C}}
\newcommand{\cA}{\mathcal{A}}
\newcommand{\cD}{\mathcal{D}}
\newcommand{\cF}{\mathcal{F}}
\newcommand{\cL}{\mathcal{L}}

\newcommand{\LF}{\mathrm{LF}}
\newcommand{\SLE}{\mathrm{SLE}}

\newcommand{\mo}[1]{{\color{purple}{#1}}}

\newtheorem{theorem}{Theorem}[section]
\newtheorem{lemma}[theorem]{Lemma}
\newtheorem{proposition}[theorem]{Proposition}
\newtheorem*{proposition*}{Proposition}
\newtheorem{corollary}[theorem]{Corollary}
\newtheorem*{corollary*}{Corollary}
\newtheorem{remark}[theorem]{Remark}
\newtheorem{prob}[theorem]{Problem}

\newtheorem{definition}[theorem]{Definition}
\newtheorem*{definitions*}{Definitions}

\newtheorem*{example*}{\bf Example}

\numberwithin{equation}{section}

\title{Reversibility of whole-plane SLE for $\kappa > 8$}
\author[1]{Morris Ang\thanks{moang@ucsd.edu}}
\author[2]{Pu Yu\thanks{py628@nyu.edu}}
\affil[1]{Department of Mathematics, University of California San Diego}
\affil[2]{Courant Institute of Mathematical Sciences, New York University}

\date{\today}
\begin{document}

\maketitle

\begin{abstract}
Whole-plane Schramm-Loewner evolution (SLE$_\kappa$) is a random fractal curve between two points on the Riemann sphere. Zhan established for $\kappa \leq 4$ that whole-plane SLE$_\kappa$ is \emph{reversible}, meaning invariant in law under conformal automorphisms swapping its endpoints. Miller and Sheffield extended this to $\kappa \leq 8$. We prove whole-plane SLE$_\kappa$ is reversible for $\kappa > 8$,  resolving the final case and answering a conjecture of Viklund and Wang.
Our argument depends on a novel mating-of-trees theorem of independent interest, where Liouville quantum gravity on the disk is decorated by an independent radial space-filling SLE curve. 
\end{abstract}



\section{Introduction}
In the past two decades, 
Schramm-Loewner evolution (SLE) has emerged as a central object of study in probability theory. SLE is a random fractal curve in the plane \cite{Sch00, RS05} describing the scaling limits of many statistical physics models at criticality \cite{smirnov-cardy, lsw-lerw-ust, ss-dgff, smirnov-ising}. It has a parameter $\kappa > 0$: when  {$\kappa \in (0,4]$} SLE is a simple curve, when $\kappa \in (4,8)$ SLE is self-intersecting but not self-crossing, and when $\kappa \geq 8$ SLE is space-filling. See for instance \cite{Law08, bn-sle-notes} for expository works on SLE.

For context, we first discuss chordal SLE, a random curve in a simply connected domain $D \subset \bbC$ from a boundary point $x$ to another boundary point $y$. We say a random curve from $x$ to $y$ is \emph{reversible} if it is invariant in law under conformal automorphisms of $D$ switching $x$ and $y$. More precisely, fixing  such a conformal automorphism $f$,  if $\eta$ is a curve from $x$ to $y$ and $\widetilde \eta$ is the time-reversal of $f \circ \eta$, then reversibility means $\eta$ and $\widetilde \eta$ agree in law up to monotone reparametrization of time.

 {The problem of SLE reversibility dates back to the very foundation of the subject. 
Schramm's definition of SLE was entirely motivated by the study of scaling limits of lattice models at criticality \cite{Sch00}: assuming a domain Markov property inherited from discrete models and the ansatz of conformal invariance, he deduced a stochastic differential equation encoding the growth of SLE. Inherent in his definition is a time-asymmetry where the starting and ending points of the curve are not interchangeable. On the other hand, many lattice models expected to converge to chordal SLE satisfy endpoint symmetry. In this way, the question of reversibility reflects a fundamental tension between the construction of SLE and its initial motivation.}

The conjecture that chordal SLE is reversible for $\kappa \in (0,8]$ was first recorded in \cite{RS05}; at the time of that conjecture, reversibility was already known for  $\kappa \in \{2, 8/3, 6, 8\}$ via scaling limits of lattice models. 
Reversibility of chordal SLE was proved by Zhan for $\kappa \in (0,4]$ \cite{zhan-chordal} and  by Miller and Sheffield for $\kappa \in (4,8)$ \cite{ig3}. On the other hand, for $\kappa > 8$ chordal SLE is not reversible \cite{RS05,   zhan2008duality}.

We now turn to 
whole-plane $\SLE_\kappa$, a random curve in $\hat \bbC := \bbC \cup \{\infty\}$ from $0$ to $\infty$. A random curve from $0$ to $\infty$ is \emph{reversible} if it is invariant in law under conformal automorphisms of $\hat \bbC$ switching $0$ and $\infty$. 
Zhan proved that whole-plane $\SLE_\kappa$ is reversible for $\kappa \leq 4$ \cite{zhan-whole-plane}, and Miller and Sheffield proved reversibility for $\kappa \in (4,8]$ \cite{MS17}. {We resolve the final case of  $\kappa>8$}.

\begin{theorem}\label{thm-main}
    Whole-plane $\SLE_\kappa$ is reversible when $\kappa > 8$. 
\end{theorem}

Theorem~\ref{thm-main} is surprising not only because of non-reversibility of chordal $\SLE_\kappa$ for $\kappa > 8$ and non-reversibility of a variant called whole-plane $\SLE_\kappa(\rho)$ for $\kappa > 8$ and $\rho > \frac\kappa2 - 4$ \cite[Remark 1.21]{MS17}, but also because it reveals a fundamental property of SLE not apparent through the lens of \emph{imaginary geometry}.  The imaginary geometry framework \cite{MS16a,MS16b,MS17,MS19} introduced by Miller and Sheffield studies SLE by coupling it with a Gaussian free field, and has proven an essential tool with wide-ranging applications such as \cite{GMS18, DMS14, KMS-removability}. The reversibility of chordal and whole-plane $\SLE_\kappa$ for $\kappa \leq 8$ can be shown by imaginary geometry \cite{ig3, MS17} (in fact, for $\kappa \in (4,8)$, this is the only known approach). However, the reversibility of whole-plane $\SLE_\kappa$ for $\kappa>8$ seems unnatural from the perspective of imaginary geometry  since the left and right boundaries of the curve interact in a complicated way  {\cite[Remark 1.22]{MS17}.}  {The reversibility for whole-plane $\SLE_\kappa(\rho)$ with $\kappa>8$ and $\rho\in(-2,\frac{\kappa}{2}-4]\backslash\{0\}$ remains an open problem.}

To our knowledge, 
apart from the illuminating work of Viklund and Wang  \cite{vw-loewner-kufarev}, there had been no reason to expect the reversibility of whole plane $\SLE_\kappa$ for $\kappa > 8$.
 They proved the inversion invariance of the $\kappa \to \infty$ large deviation rate function of whole-plane $\SLE_\kappa$, and consequently conjectured the reversibility of whole-plane $\SLE_\kappa$ for large $\kappa$. Theorem~\ref{thm-main} confirms their conjecture. 

Our arguments are substantially different from those of Zhan for $\kappa \leq 4$, who applied commutation relations for SLE \cite{zhan-chordal, zhan-whole-plane}, and Miller and Sheffield for $\kappa \leq 8$, who used imaginary geometry.  Rather, we employ the \emph{mating-of-trees} approach \cite{DMS14} where a random planar surface called \emph{Liouville quantum gravity (LQG)} is coupled with an independent SLE curve. All previously known mating-of-trees theorems \cite{DMS14, MS19, AG21} involved either chordal SLE or an SLE loop in $\hat \bbC$ or $\bbD$. 
We establish a mating-of-trees theorem for LQG on the disk coupled with radial SLE, and  for LQG on $\hat \bbC$ coupled with whole-plane SLE, resolving another conjecture of \cite{vw-loewner-kufarev}. These novel mating-of-trees theorems are noteworthy in their own right;  see for instance the survey \cite{ghs-mating-survey} for some applications of the mating-of-trees framework. 

The starting point of the original mating-of-trees theorem is the \emph{quantum zipper} coupling of reverse SLE with a certain LQG surface, from which ``zooming in'' on the base of the curve gives in the limit a forward SLE trace on a \emph{scale-invariant} LQG surface \cite{DMS14}. 
All subsequent mating-of-trees theorems were derived from the original by limiting arguments. 
 However, our radial setting is not scale-invariant, nor can it be derived from a scale-invariant picture. Our proof depends on two crucial insights. Firstly, as shown by the first author \cite{Ang23}, the quantum zipper describes dynamics on LQG surfaces arising in \emph{Liouville conformal field theory (LCFT)} \cite{DKRV16,HRV-disk}. The LCFT perspective allows us to use the quantum zipper without zooming in on a boundary point, giving us access to non-scale-invariant LQG surfaces. See \cite{AHS21, ARS21, AS21, ars-annulus, ARSZ23, ASY22} for other works that explore the interplay between LCFT and SLE.    Secondly, to pass from reverse SLE to forward SLE, we work with the infinite measure $\int_0^\infty \mathrm{raSLE}_\kappa^t \, dt$ corresponding to ``radial SLE run until a Lebesgue-typical capacity time''. This allows us to exploit the fixed-time symmetry of forward and reverse radial SLE without fixing a capacity time, which is important since capacity time is unnatural for the quantum zipper.

To prove Theorem~\ref{thm-main}, we first derive a radial mating-of-trees theorem (Theorem~\ref{thm-mot-disk-finite}) by building on the LCFT dynamics of \cite{Ang23}. Next, using a limiting argument pinching a disk into a sphere, we obtain a whole-plane mating-of-trees theorem (Theorem~\ref{thm-sphere-mot}) identifying a two-pointed  LQG sphere decorated by an independent whole-plane SLE curve with a 2D Brownian excursion. By the time-reversal symmetry of Brownian motion, 
the decorated quantum surface is invariant in law when the two points are interchanged and the curve is reversed. We conclude that whole-plane SLE is reversible.  {See Figure~\ref{fig:reversibility} for a proof summary.} Our use of mating-of-trees to prove SLE reversibility is parallel to the arguments of \cite{vw-loewner-kufarev} where  a ``mating-of-trees energy duality'' is used to establish inversion invariance of the SLE large deviation functional  {as $\kappa$ tends to infinity}. 


\medskip \noindent \textbf{Outline.}
Section~\ref{sec-prelim} gives preliminary background on LQG, Liouville conformal field theory, SLE, and mating-of-trees. In Section~\ref{sec-radial} we prove a radial mating-of-trees result (Theorem~\ref{thm-mot-disk-finite}). In Section~\ref{sec-main-proof} we take a limit to obtain a whole-plane mating-of-trees (Theorem~\ref{thm-sphere-mot}), then use it to prove Theorem~\ref{thm-main}. We  {mention related results in the literature and} list some open questions in Section~\ref{sec-open}.

\medskip 

\noindent \textbf{Acknowledgements.} 
We thank Greg Lawler, Scott Sheffield, Xin Sun, Yilin Wang and Dapeng Zhan for helpful discussions, and thank two anonymous referees for their valuable feedback.
M.A.\ was partially supported by the Simons Foundation as a Junior Fellow at the Simons Society of Fellows, and a start-up grant from the University of California San Diego. P.Y.\ was partially supported by NSF grant DMS-1712862. P.Y. thanks IAS for hosting his visit during Fall 2022.

\section{Preliminaries}\label{sec-prelim}
In this paper we work with non-probability measures and extend the terminology of ordinary probability to this setting. For a finite or $\sigma$-finite  measure space $(\Omega, \mathcal{F}, M)$, we say $X$ is a random variable if $X$ is an $\mathcal{F}$-measurable function with its \textit{law} defined via the push-forward measure $M_X=X_*M$. In this case, we say $X$ is \textit{sampled} from $M_X$ and write $M_X[f]$ for $\int f(x)M_X(dx)$. \textit{Weighting} the law of $X$ by $f(X)$ corresponds to working with the measure $d\widetilde{M}_X$ with Radon-Nikodym derivative $\frac{d\widetilde{M}_X}{dM_X} = f$. \textit{Conditioning} on some event $E\in\mathcal{F}$ (with $0<M[E]<\infty$) refers to the probability measure $\frac{M[E\cap \cdot]}{M[E]} $  on the  measurable space $(E, \mathcal{F}_E)$ with $\mathcal{F}_E = \{A\cap E: A\in\mathcal{F}\}$,  {while \emph{restricting} to $E$ refers to the measure $M[E\cap\cdot]$. }

\subsection{The Gaussian Free Field and Liouville quantum gravity}\label{subsec:def-gff}
Let $m_\bbD$ (resp.\ $m_\bbH$) be the uniform measure on the unit circle $\partial\bbD$ (resp.\ half circle $\bbH\cap\partial \bbD$).  For $X\in\{\bbD,\bbH\}$, define the Dirichlet inner product $\langle f,g\rangle_\nabla = (2\pi)^{-1}\int_X \nabla f\cdot\nabla g $ on the space $\{f\in C^\infty(X):\int_X|\nabla f|^2<\infty; \  \int f(z)m_X(dz)=0\},$ and let $H(X)$ be the closure of this space w.r.t.\ the inner product $\langle f,g\rangle_\nabla$. Let $(f_n)_{n\ge1}$ be an orthonormal basis of $H(X)$, and $(\alpha_n)_{n\ge1}$ be a collection of independent standard Gaussian variables. Then the summation
$$h_X = \sum_{n=1}^\infty \alpha_nf_n$$
a.s.\ converges in the space of distributions on $X$, and $h_X$ is the \emph{Gaussian free field (GFF)} on $X$ normalized such that $\int h_X(z)m_X(dz) = 0$. We denote its law by  $P_X$. See~\cite[Section 4.1.4]{DMS14} for more details.

Let $|z|_+ = \max\{|z|,1\}$. For $z,w\in\bar{\bbH}$, we define
$$
    G_\bbH(z,w)=-\log|z-w|-\log|z-\bar{w}|+2\log|z|_++2\log|w|_+;  \quad G_\bbH(z,\infty) = 2\log|z|_+.
$$
Similarly, for $z,w\in\bar{\bbD}$, set 
$$
G_\bbD(z,w) = -\log|z-w| - \log|1-z\bar{w}|. 
$$
Then the GFF $h_X$ is the centered Gaussian field on $X$ with covariance structure $\bbE [h_X(z)h_X(w)] = G_X(z,w)$. 

Now let $\gamma\in(0,2)$ and $Q=\frac{2}{\gamma}+\frac{\gamma}{2}$. For a conformal map $g:D\to\widetilde D$ and a generalized function $h$ on $D$, define the generalized function $g\bullet_\gamma h$ on $\widetilde{D}$ by setting {
\eqb\label{eq-coord-change}
g\bullet_\gamma h:=h\circ g^{-1}+Q\log|(g^{-1})'|
\eqe
}
A quantum surface is a $\sim_\gamma$-equivalence class of pairs $(D, h)$ where $(D,h)\sim_\gamma(\widetilde D,\widetilde h)$ if there is a conformal map $g$ with $\widetilde h = g\bullet_\gamma h$. {We call a representative $(D, h)$ an \emph{embedding} of the quantum surface. 
We will also consider quantum surfaces decorated by points and a curve; in this case we say $(D, h, \eta, (z_i)) \sim_\gamma (\widetilde D, \widetilde h, \widetilde \eta, (\widetilde z_i))$ if there is a conformal map $g: D \to \widetilde D$ such that $g \bullet_\gamma h = \widetilde h$, $g \circ \eta = \widetilde \eta$, and $g(z_i) = \widetilde z_i$ for all $i$. As before we call a representative $(D, h, \eta, (z_i))$ an embedding of the decorated quantum surface.}

For a $\gamma$-quantum surface $(D, h)$, its \textit{quantum area measure} $\mathcal{A}_h(dz)$ is defined by taking the weak limit as $\epsilon\to 0$ of $\mathcal{A}_{h_\epsilon}(dz):=\epsilon^{\frac{\gamma^2}{2}}e^{\gamma h_\epsilon(z)}dz$, where $h_\epsilon(z)$ is the circle average of $h$ over $\partial B(z, \epsilon)$. When $D=\mathbb{H}$, we can also define the  \textit{quantum boundary length measure} $\mathcal{L}_h(dx):=\lim_{\epsilon\to 0}\epsilon^{\frac{\gamma^2}{4}}e^{\frac{\gamma}{2} h_\epsilon(x)}dx$ where $h_\epsilon (x)$ is the average of $h$ over the semicircle $\{x+\epsilon e^{i\theta}:\theta\in(0,\pi)\}$. It has been shown in \cite{DS11, SW16} that all these weak limits are well-defined  for the GFF and its variants we are considering in this paper, and if $f$ is a conformal automorphism of $\bbH$ then $f_* \cA_h = \cA_{f\bullet_\gamma h}$ and $f_* \cL_h = \cL_{f\bullet_\gamma h}$. This latter point allows us to define $\mathcal{A}_h$ and $\mathcal{L}_{h}$ on other domains by conformally mapping to $\bbH$.

\subsection{The Liouville field}

Recall that $P_\bbD$ (resp.\ $P_\bbH$) is the law of the free boundary GFF on $\bbD$ (resp.\ $\bbH$) normalized to have average zero on $\partial \bbD$ (resp.\ $\partial\bbD\cap\bbH$). In the following definitions we use the shorthand $|z|_+ = \max \{ |z|, 1\}$ for $z\in \bbC$.

\begin{definition}
Let $(h, \mathbf c)$ be sampled from $P_\bbD \times [e^{-Qc} dc]$ and $\phi = h + \mathbf c$. We call $\phi$ the Liouville field on $\bbD$, and we write
$\LF_\bbD$ for the law of $\phi$.
\end{definition}

\begin{definition}
Let $(h, \mathbf c)$ be sampled from $P_\bbH \times [e^{-Qc} dc]$ and $\phi = h-2Q\log|z|_+ + \mathbf c$. We call $\phi$ the Liouville field on $\bbH$, and we write
$\LF_\bbH$ for the law of $\phi$.
\end{definition}

\begin{definition}\label{def:lf-insertion}
Let  $(\alpha,w)\in\bbR\times\bbH$ and $(\beta,s)\in\bbR\times\partial\bbH$. Let
\begin{equation*}
\begin{split}
&C_{\bbH}^{(\alpha,w),(\beta,s)} = (2\,\mathrm{Im}\,w)^{-\frac{\alpha^2}{2}}|w|_+^{-2\alpha(Q-\alpha)}|s|_+^{-\beta(Q-\frac{\beta}{2})}e^{\frac{\alpha\beta}{2}G_\bbH(w,s)}.   
\end{split}
\end{equation*}
Let $(h,\mathbf{c})$ be sampled from $C_{\bbH}^{(\alpha,w),(\beta,s) }P_\bbH\times [e^{(\alpha+\frac{\beta}{2}-Q)c}dc]$, and 
\begin{equation*}
    \phi(z) = h(z)-2Q\log|z|_++\alpha G_\bbH(z,w)+\frac{\beta}{2}G_\bbH(z,s)+\mathbf{c}.
\end{equation*}
We write $\LF_{\bbH}^{(\alpha,w),(\beta,s) }$ for the law of $\phi$ and call a sample from $\LF_{\bbH}^{(\alpha,w),(\beta,s) }$ the Liouville field on $\bbH$ with insertions ${(\alpha,w),(\beta,s) }$.
\end{definition}


\begin{definition}\label{def:LFD}
Let $\alpha,\alpha_1,\beta\in\bbR$, $w\in\bbD$ and $s\in\partial \bbD$. Let 
$$C_\bbD^{(\alpha, 0), (\alpha_1,w), (\beta, s)} = (1-|w|^2)^{-\frac{\alpha_1^2}{2}}e^{\alpha_1\alpha G_\bbD(0,w)+\frac{\alpha_1\beta}{2}G_\bbD(s,w)}.$$
Let $(h,\mathbf{c})$ be sampled from $C_\bbD^{(\alpha, 0), (\alpha_1,w), (\beta, s)}P_\bbD\times [e^{(\alpha+\alpha_1+\frac{\beta}{2}-Q)c}dc] $ and 
$$\phi(z) = h(z)+\alpha G_\bbD(z,0)+\alpha_1G_\bbD(z,w)+\frac{\beta}{2}G_\bbD(z,s)+\mathbf c.$$ We call $\phi$ the Liouville field on $\bbD$ with insertions $(\alpha,0),(\alpha_1,w),(\beta,s)$ and  write $\LF_\bbD^{(\alpha, 0), (\alpha_1,w),(\beta, s)}$ for the law of $\phi$.
\end{definition}
As we will see later in Lemma~\ref{lem:lfcoord}, the Liouville fields introduced for $\bbH$ and $\bbD$ agree up to conformal coordinate change. 

We now state the conformal covariance in $\bbH$. 
For a conformal map $f:D\to\widetilde{D}$ and a measure $M$ on $H^{-1}(D)$, let $f_*M$ be the pushforward of $M$ under the  {LQG coordinate change map $\phi\mapsto f \bullet_\gamma \phi$}. 
For $\alpha\in\bbR$, we set $\Delta_\alpha = \frac\alpha2(Q-\frac{\alpha}{2})$.
\begin{lemma}\label{lem:LF-covariance}
    Let $(\alpha,w)\in\bbR\times\bbH$ and $(\beta,s)\in\bbR\times\bbR$. Suppose $f:\bbH\to\bbH$ is a conformal map, such that $f(s)\neq\infty$. Then 
    $$\LF_{\bbH}^{(\alpha,f(w)),(\beta,f(s))} = |f'(w)|^{-2\Delta_{\alpha}}|f'(s)|^{-\Delta_{\beta}}f_*\LF_{\bbH}^{(\alpha,w),(\beta,s)}.  $$
    In particular, when $f(s)=s=0$, $f(w) = i$, we have
    \begin{equation}
        \LF_{\bbH}^{(\alpha,i),(\beta,0)} = (\mathrm{Im}\, w)^{2\Delta_{\alpha}-\Delta_\beta}|w|^{2\Delta_\beta} f_*\LF_{\bbH}^{(\alpha,w),(\beta,s)}.
    \end{equation}
\end{lemma}
\begin{proof}
This statement is proved in~\cite[Theorem 3.5]{HRV-disk}; see~\cite[Lemma 2.4]{ARS21} for an explanation.
\end{proof}

Now, we define the LCFT measure $\LF_{\bbD,\ell}^{(\alpha,0),
    (\beta,1)}$ having fixed boundary length $\ell$. 

\begin{definition}
    \label{def:disk-boundary-length}
    Let $\alpha\in\bbR, 
    \beta<Q$. Let $h$ be a sample from $P_\bbD$ and set
    $$\widetilde h(z) = h +\alpha G_\bbD(z,0)+
    \frac{\beta}2G_\bbD(z,1).$$ 
    Fix $\ell>0$, and let $L=\mathcal{L}_{\widetilde h}(\partial\bbD)$. Define the measure $\LF_{\bbD,\ell}^{(\alpha,0),
    (\beta,1)}$ to be the law of $\widetilde h+\frac{2}{\gamma}\log\frac{\ell}{L}$ under the reweighted measure $\frac{2}{\gamma}\frac{\ell^{\frac{2\alpha+\beta-2Q}{\gamma}-1}}{L^{\frac{2\alpha+\beta-2Q}{\gamma}}}P_\bbD(dh)$.
\end{definition}

 \begin{lemma}\label{lem:disk-boundary-length}
     In the setting of Definition~\ref{def:disk-boundary-length}, {$\{\LF_{\bbD, \ell}^{(\alpha,0),(\beta,1)}\}_{\ell > 0}$ is a \emph{disintegration} of $\LF_\bbD^{(\alpha,0),(\beta,1)}$ over its boundary length.} That is, any sample $\phi$ from $\LF_{\bbD,\ell}^{(\alpha,0),(\beta,1)}$ has $\mathcal{L}_\phi(\partial\bbD)=\ell$, and 
     \begin{equation}\label{eq:lf-disintegrate}
         \LF_\bbD^{(\alpha,0),
         (\beta,1)} = \int_0^\infty \LF_{\bbD,\ell}^{(\alpha,0),
         (\beta,1)}d\ell.
     \end{equation}
     Moreover, if $\alpha
     +\frac{\beta}{2}>Q$, we have $|\LF_{\bbD,\ell}^{(\alpha,0),
     (\beta,1)}| = C\ell^{\frac{2\alpha+
     \beta-2Q}{\gamma}-1}$ for some finite constant $C$.
 \end{lemma}
\begin{proof}
    {First, $\cL_\phi(\partial \bbD) = \cL_{\widetilde h + \frac2\gamma \log \frac \ell L}(\partial \bbD) = \frac \ell L \cL_{\widetilde h }(\partial \bbD) = \ell$.} Next, for any nonnegative measurable function $F$ on $H^{-1}(\mathbb{D})$,
    $$\int_0^\infty\int F(\widetilde{h}+\frac{2}{\gamma}\log\frac{\ell}{L})\frac{2}{\gamma}\frac{\ell^{\frac{2\alpha+
    \beta-2Q}{\gamma}-1}}{L^{\frac{2\alpha+
    \beta-2Q}{\gamma}}}P_{\mathbb{D}}(dh)d\ell = \int_{\mathbb{R}}\int F(\widetilde{h}+c)e^{(\alpha+\frac{\beta}{2}-Q)c}P_\bbD(dh) dc$$
    using Fubini's theorem and the change of variables $c = \frac{2}{\gamma}\log\frac{\ell}{L}$. This justifies~\eqref{eq:lf-disintegrate}. For the last claim,
    \begin{equation}
    \begin{split}
    \LF_\bbD^{(\alpha,0),
    (\beta,1)}[\{{\mathcal{L}_\phi(\partial\bbD)\in[a,b]}\}] &=  \int\int \mathds{1}_{e^{\frac{\gamma}{2}c}L\in[a,b]}e^{(\alpha+
    \frac{\beta}{2}-Q)c} P_\bbD(dh)dc\\&= \frac{2}{\gamma}\int L^{-\frac{2\alpha+
    \beta-2Q}{\gamma}} P_\bbD(dh)\cdot \int_a^b\ell^{\frac{2\alpha+
    \beta-2Q}{\gamma}-1}d\ell
    \end{split}\end{equation}
    where we used the change of variables $\ell = e^{\frac{\gamma}{2}c}L$. Since $\alpha+
    \frac{\beta}{2}>Q$, the integral $\int L^{-\frac{2\alpha+
    \beta-2Q}{\gamma}} P_\bbD(dh)$ is finite (see e.g.~\cite{HRV-disk,RZ20b}) and the claim then follows.
\end{proof}

{
As we see next, sampling a point from the LQG area measure corresponds to adding an LCFT insertion of size $\gamma$.} Recall $\mathcal{A}_\phi(dz)$ denotes the quantum area measure.
\begin{lemma}\label{lem:qtypical}
   Let $w\in\bbD$, $\alpha,\beta\in\bbR$ and $s\in\partial\bbD$. Then we have
   $$\mathcal{A}_\phi(dz)\LF_\bbD^{(\alpha,0),(\beta,s)}(d\phi) =  \LF_\bbD^{(\alpha,0),(\beta,s),(\gamma,z)}(d\phi)dz  .$$
\end{lemma}
{
\begin{proof}
The proof is identical to that of~\cite[Proposition 2.5]{Ang23}.
\end{proof}
} 
Finally, the Liouville fields on $\bbH$ and $\bbD$ agree up to coordinate change; we now verify the case that we need for this paper.
\begin{lemma}\label{lem:lfcoord}
     Let $\alpha,\beta\in\bbR$ with $\alpha+\frac{\beta}{2}=Q$. For $w\in\bbH$, let $g:\bbH\to\bbD$ be a conformal map with $g(w)=0$ and $g(0)=1$. Then 
    \begin{equation}\label{eq:lem-lfcoord}
        \LF_{\bbD}^{(\alpha,0),(\beta,1)} = 2^{\frac{\alpha^2}{2}}(\mathrm{Im}\, w)^{2\Delta_{\alpha}-\Delta_\beta}|w|^{2\Delta_\beta} g_*\LF_{\bbH}^{(\alpha,w),(\beta,0)}.
    \end{equation}
\end{lemma}
\begin{proof}
    We will show the claim for $w = i$, then the general case follows by using Lemma~\ref{lem:LF-covariance}.

    Let $g: \bbH \to \bbD$ be the conformal map such that $g(i) = 0, g(0) = 1$. Explicitly, it is given by $g(z) = \frac{i - z}{i+z}$.
    By the conformal invariance of the free boundary GFF viewed as a distribution modulo additive constant, if 
    $(h_\bbH, \mathbf c_\bbH) \sim P_\bbH \times dc$ and $(h_\bbD, \mathbf c_\bbD) \sim P_\bbD \times dc$, then $h_\bbH + \mathbf c_\bbH \stackrel d= (h_\bbD + \mathbf c_\bbD) \circ g$. 
    Next, 
    using the formulas for $G_\bbH$ and $G_\bbD$ in Section~\ref{subsec:def-gff}, one can directly check that for some constant $C$, we have 
    \[\alpha G_\bbH(\cdot, i) + \frac\beta2 G_\bbH(\cdot, 0) -2Q \log |\cdot|_+ = (\alpha G_\bbD(\cdot, 0) + \frac\beta2 G_\bbD (\cdot, 1) )\circ g + Q \log |g'| + C \quad \text{ for all }z \in \bbH.\]
    Combining this with the translation invariance of Lebesgue measure, we conclude that 
    \[g \bullet_\gamma (h_\bbH + \mathbf c_\bbH + \alpha G_\bbH(\cdot, i) + \frac\beta2 G_\bbH(\cdot, 0) -2Q \log |\cdot|_+ ) \stackrel d= h_\bbD + \mathbf c_\bbD + \alpha G_\bbD(\cdot, 0) + \frac\beta2 G_\bbD (\cdot, 1). \]
    Thus~\eqref{eq:lem-lfcoord} holds for $w = i$, as needed. 
\end{proof}


\subsection{Forward and reverse SLE}

In this section we briefly recall the forward and reverse radial SLE processes, and whole-plane SLE. We will not give precise definitions since they will not be used later, but curious readers can refer to \cite{Law08}.


Forward radial $\SLE_\kappa$ in $\bbD$ from $1$ to $0$ is a random non-self-crosing curve $\eta:[0,\infty) \to \ol \bbD$ with $\eta(0) = 1$ and $\lim_{t \to \infty} \eta(t) = 0$. Let $K_t$ be the compact subset of $\ol \bbD$ such that $\ol \bbD \backslash K_t$ is the connected component of $\bbD \backslash \eta([0,t])$ containing $0$, and let $g_t: \bbD \backslash K_t \to \bbD$ be the conformal map with $g_t(0) = 0$ and $g_t'(0) > 0$. The curve $\eta$ is parametrized by log conformal radius, meaning that for each $t$ we have $g_t'(0) = e^t$. It turns out that there is a random process $U_t \stackrel d=e^{i\sqrt{\kappa}B_t}$ (where $B_t$ is standard Brownian motion) such that 
\eqb\label{eq-rad-sle}
dg_t(z) = \Phi(U_t, g_t(z))\,dt\quad \text{ for } z\in\mathbb{D}\backslash K_t  \text{ and }\Phi(u,z) := z\frac{u+z}{u-z}.
\eqe
{In fact,~\eqref{eq-rad-sle} 
and the initial condition $g_0(z) = z$ define the family  of conformal maps $(g_t)_{t \geq 0}$ and hence radial $\SLE_\kappa$}, see \cite{Law08} for details.

Similarly, whole-plane $\SLE_\kappa$ is a random non-self-crossing curve $\eta:(-\infty, \infty) \to \bbC$ {from $0$ to $\infty$}, such that if $K_t$ is the compact set such that $\bbC \backslash K_t$ is the unbounded connected component of $\bbC \backslash \eta((-\infty, t))$, and $g_t: \bbC \backslash K_t \to \bbC \backslash \bbD$ is the conformal map such that $g_t(\infty) = \infty$ and $g_t'(\infty) > 0$, then $$dg_t(z) = \Phi(U_t, g_t(z))\,dt \quad \text{ for } z \in \bbC \backslash K_t$$ where $U_t\stackrel d= e^{i\sqrt{\kappa}B_t}$ and $(B_t)_{t\in\bbR}$ is two-sided standard Brownian motion. This curve extends continuously to its starting and ending points, i.e.  $\lim_{t \to -\infty} \eta(t) = 0$ and $\lim_{t \to \infty} \eta(t) = \infty$~\cite{Lawler11continuity,MS17}.

Now we discuss \emph{centered reverse} radial SLE. Unlike the forward case where we have a single random curve, centered reverse radial SLE is a random process of curves $(\eta_t)_{t \geq 0}$. Each curve $\eta_t:[0,t] \to \ol \bbD$ is parametrized by log conformal radius and has starting point $\eta_t(0) = 1$, and $(\eta_t)_{t \geq 0}$ satisfies the compatibility relation that for $s < t$, if $\widetilde f_{s,t}$ is the conformal map from $\bbD$ to the connected component of {$\bbD \backslash \eta_t([0,t-s])$} containing 0  such that $\widetilde f_{s,t}(1) = \eta_t(t-s)$ and $\widetilde f_{s,t}(0) = 0$, then $\eta_s = \widetilde f_{s,t}^{-1} \circ \eta_t(\cdot + t-s)|_{[0,s]}$. Informally, this compatibility relation means that the process $(\eta_t)_{t \geq 0}$ grows from the base of the curve. 
We call $\widetilde f_{0, t}$ the centered reverse Loewner map.  
The process $(\eta_t)_{t \geq 0}$ satisfies the stochastic differential equation 
\eqb\label{eq-centered-reverse-radial}
d\widetilde f_{0,t}(z) = -i \sqrt\kappa \widetilde f_{0,t}(z) dB_t -\Phi(1, \widetilde f_{0,t}(z))\,dt\ \text{for}\ z\in\ol{\mathbb{D}}.
\eqe
One can show via the time-reversal symmetry of Brownian motion that for each fixed $t$, the curve $\eta_t$ has the law of forward radial SLE run for time $t$. 

For $z_0 \in \bbH$ and $\rho \in \bbR$, there is also a random process $(\eta_t)_{t \geq 0}$ called \emph{centered reverse chordal $\SLE_\kappa(\rho)$ with force point at $z_0$} (see e.g.~\cite[Section 4.3]{zhanrohde2016}, ~\cite[Section 3.3.1]{DMS14}). Each $\eta_t:[0,t] \to \ol \bbH$ is parametrized by half-plane capacity, has $\eta_t(0) = 0$, and satisfies a compatibility relation analogous to that of the radial case. It is defined by a stochastic differential equation similar to~\eqref{eq-centered-reverse-radial} which we omit here. For each $t>0$ let $\widetilde f_{\bbH,t}: \bbH \to \bbH \backslash \eta_t$ be the conformal map with $\widetilde f_{\bbH,t}(0) = \eta_t(t)$ and $\widetilde f_{\bbH,t}(z) = z + O(1)$ as $z \to \infty$; we call $\widetilde f_{\bbH,t}$ the  centered reverse Loewner map.

Finally, \cite[Theorem 4.6]{zhanrohde2016}\footnote{{They use a different notation for weights of force points, see Remark 2 immediately after \cite[Corollary 4.8]{zhanrohde2016}.}} gives a change of coordinates result for  reverse chordal $\SLE$:
\begin{lemma}\label{lem:rslecoord}
    Fix $\kappa>0$. Let $(\eta_t)_{t \geq 0}$ be a centered reverse chordal $\SLE_\kappa(\kappa + 6)$ process with force point at $\widetilde z_0\in\bbH$. Let $\widetilde f_t$ be its associated reverse centered Loewner map. Let $\varphi_0:\bbH\to\bbD$ be the conformal map with $\varphi_0(\widetilde z_0)=0$ and $\varphi_0(0) = 1$, and $\varphi_t:\bbH\to\bbD$  the conformal map such that  $\varphi_t(\widetilde f_t(\widetilde z_0)) = 0$ and $\varphi_t(0) = 1$. 
    Let $\eta_t'$ be $\varphi_t \circ \eta_t$ parametrized by log conformal radius.
    Then up to a time change, $(\eta_t')_{t\ge0}$ has the law of centered reverse  radial $\SLE_\kappa$  stopped at the time  $\varphi_0(\infty)$ hits the driving function, i.e.\ the first time $s$ when $\widetilde f_{0,s}(\varphi_0(\infty)) = 1$ where $\widetilde f_{0,s}$ is the centered  reverse Loewner map of the reverse radial $\SLE_\kappa$.
\end{lemma}


\subsection{Chordal mating-of-trees and special quantum surfaces}
In this section we state the chordal mating-of-trees theorem of \cite{DMS14}, and recall the definition of the \emph{quantum cone} from~\cite{She16a,DMS14} and the \emph{quantum cell} from~\cite{Ang23}.

Let $\cC = (\bbR \times [0,2\pi])/{\sim}$ be the horizontal cylinder obtained by gluing the upper and lower boundaries of the strip via the identification $x \sim x + 2\pi i$.
We define the GFF on $\cC$ as in Section~\ref{subsec:def-gff}, with $m_{\mathcal C}$ the uniform measure on $(\{0\} \times [0,2\pi])/{\sim}$, and likewise define the Hilbert space $H(\mathcal C)$. As explained in, e.g., \cite[Section 4.1.7]{DMS14}, we may decompose $H(\mathcal{C}) = H_{\text{av}}(\mathcal C)\oplus H_{\text{lat}}(\mathcal C)$, where $H_{\text{av}}(\mathcal C)$ (resp.\ $H_{\text{lat}}(\mathcal C)$) is the subspace of functions which are constant (resp.\ have mean 0) on $\{t\}\times [0,2\pi]$ for each $t\in\bbR$. This gives a decomposition $h_{\mathcal{C}} = h_{\text{av}}+h_{\text{lat}}$ of $h_{\mathcal C}$ into two independent components.

Now we introduce   the $\gamma$-LQG surfaces  called quantum cones  {via an embedding in $(\cC, -\infty, +\infty)$. Near $-\infty$ it has finite quantum area, but every neighborhood of $+\infty$ has infinite quantum area.}
\begin{definition}[$\alpha$-quantum cone]
    Fix $\alpha < Q$. Suppose $\psi_{\textup{av}}$ and $\psi_{\textup{lat}}$ are independent distributions on $\mathcal{S}$ such that: 
	\begin{itemize}
		\item We have $\psi_{\textup{av}}(z)= X_{{\mathrm{Re}}\, z}$ for $z \in \mathcal C$, where
		\begin{equation}
		X_t:=\left\{ \begin{array}{rcl} 
			B_{t}-(Q-\alpha) t & \mbox{for} & t\ge 0\\
			\widetilde{B}_{-t} +(Q-\alpha) t & \mbox{for} & t<0
		\end{array} 
		\right.
		\end{equation}
        and  $(B_t)_{t\ge 0}$ and $(\widetilde{B}_t)_{t\ge 0}$ are independent standard Brownian motions conditioned on  $ \widetilde{B}_{t} - (Q-\alpha)t<0$ for all $t>0$\footnote{This conditioning can be made sense via Bessel processes; see e.g.~\cite[Section 4.2]{DMS14}.};  
		\item $\psi_{\textup{lat}}$ has the same law as $h_{\textup{lat}}$.
  \end{itemize}
	 Set $\psi=\psi_{\textup{av}}+\psi_{\textup{lat}}$. We call $(\mathcal{C}, \psi, -\infty, +\infty)/{\sim_\gamma}$ an \emph{$\alpha$-quantum cone}.
\end{definition}

For $\kappa >4$, there is a random curve in $\bbC$ called \emph{space-filling $\SLE_\kappa$ from $\infty$ to $\infty$}. It is defined via the imaginary geometry flow lines of a whole-plane GFF. Space-filling $\SLE_\kappa$ from $\infty$ to $\infty$ is reversible since its construction is symmetric. Moreover, if $\kappa \geq 8$, for each $z \in \bbC$ the regions covered by the curve before and after hitting $z$ are simply connected, and conditioned on the curve up until it hits $z$, it subsequently evolves as chordal $\SLE_\kappa$ from $z$ to $\infty$ in the complementary domain.
This follows from the flow line construction of space-filling $\SLE_\kappa$, see~\cite[Section 1.2.3]{MS17} for more details.

We are ready to state the mating-of-trees theorem ~\cite[Theorem 1.9, Theorem 1.11]{DMS14}. We shall focus on the $\kappa>8$ regime.
\begin{theorem}\label{thm:mot-chordal}
    Let $\kappa>8$ and $\gamma = \frac{4}{\sqrt{\kappa}}$. Let $(\bbC,\phi,0,\infty)$ be an embedding of a  $\gamma$-quantum cone and $\eta$ an independent space-filling $\SLE_\kappa$ curve from $\infty$ to $\infty$, and we reparameterize $\eta$ by the $\gamma$-LQG measure, in the sense that $\eta(0)=0$ and $\mathcal{A}_\phi(\eta([s,t]))=t-s$ for $-\infty<s<t<\infty$. Define $X_t^-, X_t^+, Y_t^-, Y_t^+$ as in Figure~\ref{fig:quantumcell} (left, middle) and let $X_t := X_t^+ - X_t^-$ and $Y_t := Y_t^+ - Y_t^-$.  Then $(X_t, Y_t)_{t\in\bbR}$ is a correlated two-sided two-dimensional Brownian motion with $X_0=Y_0=0$, with covariance
\begin{equation}\label{eqn-thm:mot-chordal}
    \var(X_t) = \var(Y_t) = \mathbbm{a}^2|t|; \quad \mathrm{cov}(X_t,Y_t) = -\cos(\frac{4\pi}{\kappa}) \mathbbm{a}^2|t| \quad \text{where }\mathbbm{a}^2:=\frac2{\sin(\frac{4\pi}{\kappa})}.
\end{equation}
Moreover, the pair $(X, Y)$ a.s.\ determines {the decorated quantum surface $(\bbC, \phi, \eta,0, \infty)/{\sim_\gamma}$.}
\end{theorem}
We can interpret $X_t$ (resp.\ $Y_t$) as the change in the quantum length of the left (resp.\ right) boundary
of $\eta$ relative to time 0.
The covariance in~\eqref{eqn-thm:mot-chordal} was computed in~\cite{GHMS15} while the constant $\mathbbm{a}$ was obtained in~\cite{ARS21}.

Let $(\phi,\eta)$ and $(X,Y)$ be as in the statement of Theorem~\ref{thm:mot-chordal}. For each $a > 0$, let $D_a = \eta([0, a])$, $p = \eta(0)=0$ and $q = \eta(a)$. Let $x_L$
(resp.\ $x_R$) be the last point on the left (resp. right) boundary arc of $\eta((-\infty,0])$ hit by $\eta$ before time $a$. See Figure~\ref{fig:quantumcell} (right).

\begin{definition}\label{def-cell}
    We call the $\SLE_\kappa$-decorated quantum surface $\mathcal{C}_a := (D_a,h,\eta|_{[0,a]};p,q,x_L,x_R)/{\sim_\gamma}$ an area $a$ quantum cell, and denote its law by $P_a$. We call $(X_t, Y_t)_{[0,a]}$ its boundary length process,  {and $X_a^- = - \inf_{0<t<a} X_t$, $X_a^+ = X_a + X_a^-$, $Y_a^- = -\inf_{0<t<a} Y_t$, $Y_a^+ = Y_a + Y_a^-$ its boundary lengths.}
\end{definition}
Note that the quantum length of the arc between $p$ and $x_L$ (resp.\ $x_R$) is $X_a^-$ (resp.\ $Y_a^-$), and the quantum length of the arc between $q$ and $x_L$ (resp.\ $x_R$) is $X_a^+$ (resp.\ $Y_a^+$).
\cite{Ang23} gives a different but equivalent definition of the quantum cell in terms of the so-called weight $2-\frac{\gamma^2}2$ quantum wedge; the equivalence follows from the fact that in the setting of Theorem~\ref{thm:mot-chordal}, the quantum surface $(\eta((0, \infty)), \phi, 0, \infty)/{\sim}$ has the law of the weight $2-\frac{\gamma^2}2$ quantum wedge \cite[Theorem 1.9]{DMS14}. 

\begin{figure}
    \centering
    \includegraphics[scale = 0.53]{figures/quantumcell2.pdf}
    \caption{\textbf{Left:} Let $t>0$. Let $X_t^-$ be the quantum length along $\partial (\eta((-\infty, 0)))$ from $0$ to the leftmost point of $\eta([0,t]) \cap \partial (\eta((-\infty, 0)))$, and $X_t^+$ the quantum length of the counterclockwise boundary arc of $\eta([0,t])$ from this point to $\eta(t)$. Likewise define $Y_t^-, Y_t^+$. 
    \textbf{Middle:} When $t < 0$ we let $X_t^-$ be the quantum length along $\partial (\eta((-\infty, t)))$ from $\eta(t)$ to the leftmost point of $\eta([t,0]) \cap \partial (\eta((-\infty, t)))$, and $X_t^+$ the quantum length of the counterclockwise boundary arc of $\eta([t,0])$ from this point to $0$. Likewise define $Y_t^-, Y_t^+$. 
    \textbf{Right:} 
    An illustration of a quantum cell of quantum area $a$ and its boundary lengths $X_a^+,X_a^-,Y_a^+,Y_a^-$.}
    \label{fig:quantumcell}
\end{figure}

By~\cite[Remark 2.9]{Ang23},  $\mathcal{C}_a$ is measurable with respect to $(D_a,h,\eta|_{[0,a]})/{\sim_\gamma}$ since $\kappa > 8$, and therefore we will often omit the marked points of $\mathcal{C}_a$ for notational simplicity. 
The quantum surface $(D_a,h,\eta|_{[0,a]})/{\sim_\gamma}$ is measurable with respect to $(X_t,Y_t)_{0\le t\le a}$ \cite[Lemma 2.17]{AG21}, and we denote the map sending $(X_t,Y_t)_{0\le t\le a}$ to $(D_a,h,\eta|_{[0,a]})/{\sim_\gamma}$ by $F$. We now give two properties of $F$.

\begin{lemma}[Reversibility of $F$]\label{lem-F-reversible}
Fix $a>0$, sample $\mathcal C_a = (D, h, \eta)/{\sim_\gamma}$ from $P_a$, and let $(X_t, Y_t)_{[0,a]}$ be its boundary length process, so $F((X_t, Y_t)_{[0,a]}) = \mathcal C_a$ a.s.. Let $\widetilde {\mathcal C}_a = (D, h, \widetilde \eta)/{\sim_\gamma}$ where $\widetilde \eta$ is the time-reversal of $\eta$, and let $(\widetilde X_t, \widetilde Y_t)_{[0,a]} =  {(Y_{a-t}, X_{a - t})_{[0,a]}}$ be the time-reversal of $(X_t, Y_t)_{[0,a]}$. Then $F((\widetilde X_t, \widetilde Y_t)_{[0,a]}) = \widetilde {\mathcal C}_a$ a.s..
\end{lemma}
\begin{proof}
    Let $(\bbC, h, 0, \infty)$ be an embedding of a  $\gamma$-quantum cone and let $\eta$ be an independent SLE from $\infty$ to $\infty$ in $\bbC$ parametrized by quantum area such that $\eta(0) = 0$.
    Let $\mathcal C_a = (\eta([0,a], h, \eta|_{[0,a]})$ so the law of $\cC_a$ is $P_a$, and let $(X_t, Y_t)_{[0,a]}$ be its boundary length process. 
    Let $\widetilde \eta$ be the time-reversal of $\eta$, then by the reversibility of $\SLE$ from $\infty$ to $\infty$ in $\bbC$ we have $(\bbC, h, \eta, 0, \infty)/{\sim_\gamma} \stackrel d= (\bbC, h, \widetilde\eta, 0, \infty)/{\sim_\gamma}$. Let $\widetilde \eta'(\cdot) = \widetilde \eta(\cdot -a)$ (so $\widetilde \eta'|_{[0,a]}$ is the time-reversal of $\eta|_{[0,a]}$), then \cite[Lemma 8.3]{DMS14} implies $(\bbC, h, \widetilde\eta, 0, \infty)/{\sim_\gamma} \stackrel d= (\bbC, h, \widetilde\eta', \eta(a), \infty)/{\sim_\gamma}$, that is, {$(\bbC, h, \widetilde\eta', \eta(a), \infty)/{\sim_\gamma}$} is a quantum cone decorated by an independent SLE from $\infty$ to $\infty$ in $\bbC$. We conclude that the law of $\widetilde \cC_a$ is also $P_a$, and directly from the definition of boundary length process, the boundary length process of $\widetilde \cC_a$ is $(\widetilde X_t, \widetilde Y_t)_{[0,a]}$, so $F((\widetilde X_t, \widetilde Y_t)_{[0,a]}) = \widetilde \cC_a$ a.s.. 
\end{proof}

\begin{lemma}[Concatenation compatibility of $F$]\label{lem-concat}
Let $a_1, a_2 > 0$, and let  $(X_t, Y_t)_{t \in \bbR}$ be as in~\eqref{eqn-thm:mot-chordal}. 
Let $\cC_1 = F((X_t,Y_t)_{[0,a_1]})$, let $\cC_2 = F((X_{t+a_1} - X_{a_1}, Y_{t+a_1} - Y_{a_1})_{[0,a_2]})$, and let $\cC = F((X_t, Y_t)_{[0,a_1+a_2]})$. Almost surely, $\cC_1$ and $\cC_2$ are the curve-decorated quantum surfaces obtained from $\cC$ by restricting to the domains parametrized by its curve on the time intervals $[0,a_1]$ and $[a_1, a_1+a_2]$. 
\end{lemma}
\begin{proof}
This is immediate from the definition of $F$ and the fact that if $(\bbC, \phi, 0, \infty)$ is an embedding of a $\gamma$-quantum cone and $\eta$ is an independent space-filling $\SLE_\kappa$ from $\infty$ to $\infty$ parametrized by quantum area, then $(\bbC, \phi, \eta, 0, \infty)/{\sim_\gamma} \stackrel d= (\bbC, \phi, \eta(\cdot + a_1), \eta(a_1), \infty)/{\sim_\gamma}$ \cite[Lemma 8.3]{DMS14}.
\end{proof}

Finally, we recall the definition of the quantum sphere of \cite{DMS14}.  {This is a two-pointed quantum surface with finite quantum area.}
\begin{definition}\label{def-sphere}
    Let $\alpha < Q$. Let $(B_s)_{s \geq 0}$ be a standard Brownian motion conditioned on $B_s - (Q-\alpha)s < 0$ for all $s>0$, and let $(\widetilde B_s)_{s \geq 0}$ be an independent copy of $(B_s)_{s \geq 0}$. Let
    \[Y_t = \left\{
	\begin{array}{ll}
		B_t - (Q-\alpha)t  & \mbox{if } t \geq 0 \\
		\widetilde B_{-t} + (Q-\alpha)t & \mbox{if } t < 0
	\end{array}
\right.
\]
Let $h^1(z) = Y_{\mathrm{Re}\, z}$ for $z \in \cC$, and let $h^2$ be independent of $h^1$ and have the law of the lateral component of the GFF on $\cC$. Let $\hat h = h^1 + h^2$. Let $\mathbf c \in \bbR$ be independently sampled from $\frac\gamma2 e^{2(\alpha-Q)c}\,dc$. Let $\mathcal M_2^\mathrm{sph}(\alpha)$ be the infinite measure describing the law of the decorated quantum surface $(\cC, \hat h + \mathbf c, -\infty, +\infty)/{\sim_\gamma}$.
\end{definition}

 {
\subsection{LCFT and the quantum zipper}\label{sec-quantum-zipper}
In this section we state a special case of the chordal  quantum zipper for LCFT obtained in \cite{Ang23}. It will be used in Section~\ref{sec-radial-zipper} to derive a radial quantum zipper for LCFT. 

Let $\kappa > 8$ and $\gamma = \frac4{\sqrt\kappa}$. Let $\mathrm{BM}_\kappa$ denote the law of (one-sided) correlated two-dimensional Brownian motion $(X_t, Y_t)_{t \geq 0}$ with $X_0 = Y_0=0$ and covariance given by~\eqref{eqn-thm:mot-chordal}.
Let $\alpha \in \bbR$, and sample $(\tilde \psi_0, (X_t, Y_t)_{t\geq 0}) \sim \LF_\bbH^{(\alpha, i), (-\frac\gamma2, 0)} \times \mathrm{BM}_\kappa$. For each $s>0$ let $\cC_s = F((X_\cdot, Y_\cdot)|_{[0,s]})$, and on the event that $-\inf_{u< s} X_u - \inf_{u< s} Y_u < \cL_{\psi_0}(\bbR)$, conformally weld $(\bbH,  \psi_0, i, 0, \infty)/{\sim_\gamma}$ to $\cC_s$ by identifying the first marked points of each quantum surface and identifying the two boundary arcs of $ \cC_s$ adjacent to its first marked point to a boundary interval of $(\bbH, \psi_0, i, 0)/{\sim_\gamma}$; this identification is via the quantum length measures of the two quantum surfaces. The resulting curve-decorated quantum surface has many possible embeddings in $\bbH$. Let $(\bbH, \tilde {\psi}_s, \tilde{ \eta}_s)$ be the unique embedding such that if $f:  \bbH \to \bbH\backslash \tilde{\eta}_s $ is the conformal map which fixes $\infty$, sends the tip of $\tilde{\eta}_s$ to $0$, and satisfies $f(z) = z + O(1)$ as $z \to \infty$, then $\tilde \psi_0 = f^{-1} \bullet_\gamma \tilde \psi_s$. In this way, we obtain a process $(\wt \psi_s, \wt \eta_s)$. Let $(\psi_t, \eta_t)$ be the monotone reparametrization of the process such the half-plane capacity of the trace of $\eta_t$ is $2t$. 
\begin{lemma}\label{lem-chordal-zipper}
	For any stoppping time $\sigma$ for the filtration $\cF_t = \sigma(\eta_t)$, the law of $(\psi_\sigma, \eta_\sigma)$ is 
	\[
		\frac{1}{\mathcal{Z}_\alpha(\widetilde{f}_{\bbH, \sigma}(i))}\LF_\bbH^{(\alpha,\widetilde f_{\bbH, \sigma}(i)), (-\frac{\gamma}{2},0)}\,\mathrm{rSLE}_{\kappa,2\sqrt\kappa \alpha}^{\sigma}, \qquad \mathcal{Z}_\alpha(z) := (2\mathrm{Im}\,z)^{-\frac{\alpha^2}2} |z|^{\frac2{\sqrt\kappa}\alpha},
	\]
	where $\mathrm{rSLE}_{\kappa,\rho}^{\sigma}$ denotes the law of centered reverse chordal $\SLE_\kappa(\rho)$ with the force point located at $i$ run until the stopping time $\sigma$ and $\widetilde f_{\bbH,\sigma}$ is its associated reverse Loewner map.
\end{lemma}
\begin{proof}
	This is the special case of \cite[Theorem 1.8]{Ang23} where there is a single bulk insertion and a single boundary insertion, phrased in terms of centered reverse chordal $\SLE_\kappa(\rho)$ rather than reverse chordal $\SLE_\kappa(\rho)$. 
\end{proof}
}

\section{A radial mating-of-trees theorem}\label{sec-radial}

In this section, we prove our radial mating-of-trees result Theorem~\ref{thm-mot-disk-finite}. Throughout this section, let $\gamma\in(0,\sqrt{2})$ and $\kappa=\frac{16}{\gamma^2}>8$.


Sample {$\phi\sim \LF_\bbD^{(Q-\frac\gamma4, 0), (\frac{3\gamma}2, 1)}$} conditioned on having quantum boundary length 1, let $A = \cA_\phi(\bbD)$, and let $\eta:[0,A] \to \ol \bbD$ be an independent radial $\SLE_\kappa$ in $\bbD$ from 1 to 0 parametrized by its $\mathcal{A}_\phi$-quantum area.
There is a unique continuous process
$(X_t, Y_t)_{[0,A]}$ starting at $(X_0, Y_0) = (0,0)$ which keeps track of the local changes in the left and right LQG boundary lengths of $\bbD \backslash \eta([0,t])$ in the following sense. For any time $s \in (0,A)$ and any point $p \in \partial (\bbD \backslash \eta([0,s]))$ different from $\eta(s)$, let $\sigma>s$ be the next time $\eta$ hits $p$. For each  $t \in [s, \sigma)$, let $X_t^s$ (resp.\ $Y_t^s$) be the quantum length of the clockwise (resp.\ counterclockwise) boundary arc of $\bbD \backslash \eta([0,t])$ from $\eta(t)$ to $p$. Then $(X_t - X_s, Y_t - Y_s)_{[s, \sigma)} = (X_t^s - X_s^s, Y_t^s - Y_s^s)_{[s,\sigma)}$. See Figure~\ref{fig-lengths} (left, middle). This process can be constructed on the time interval $[0,A)$  by shifting the point $p$ countably many times, and its value at $A$ is defined by taking a limit. Note that these LQG lengths exist and are finite by local absolute continuity with respect to the setting of Theorem~\ref{thm:mot-chordal}.

\begin{theorem}[{$\kappa>8$} radial mating-of-trees]\label{thm-mot-disk-finite}
The process $(X_t,Y_t)_{0\le t\le A}$ has the law of 2-dimensional Brownian motion with covariance~\eqref{eqn-thm:mot-chordal} stopped at the first time that $1 + X_\cdot + Y_\cdot = 0$. Moreover, for $0\le s<t$, on the event that $t<A$ and $\eta([s,t])$ is simply connected, we have 
\eqb\label{eq-thm-radial}
F((X_{\cdot + s}-X_{s}, Y_{\cdot + s}-Y_{s})|_{[0, t-s]}) = (\eta([s, t]), \phi, \eta(\cdot + s)|_{[0,t-s]})/{\sim_\gamma}\quad  \text{almost surely}.\eqe
Here, $F$ is as in Lemma~\ref{lem-F-reversible}.
\end{theorem}

{We note that when $\eta([s,t])$ is not simply connected, then instead the right hand side of~\eqref{eq-thm-radial} is obtained from the left hand side by conformally welding its boundary to itself.}
In Theorem~\ref{thm-mot-disk-finite} the curve-decorated quantum surface $(\bbD, \phi, \eta, 0, 1)/{\sim_\gamma}$ can a.s.\ be recovered from $(X_t, Y_t)_{[0,A]}$ by conformally welding countably many {simply connected} quantum surfaces of the form $(\eta([s,t]), \phi, \eta(\cdot + s)|_{[0, t-s]})/{\sim_\gamma}$, each of which is measurable with respect to $(X_t, Y_t)_{0\leq t \leq A}$ by~\eqref{eq-thm-radial}.

 
\begin{figure}
    \centering
   \includegraphics[scale=0.4]{figures/lengths.pdf}
    \caption{\textbf{Left:} The boundary length process $(X_t, Y_t)_{[0,A]}$ of Theorem~\ref{thm-mot-disk-finite} is characterized by $X_0 = Y_0 = 0$ and the property that for each time $s$ and choice of boundary point $p \in \partial (\bbD \backslash \eta([0,s]))$ different from $\eta(s)$, for any time $t > s$ before  the time $\eta$ next hits $p$, we have $(X_t^s - X_s^s, Y_t^s - Y_s^s) = (X_t - X_s, Y_t - Y_s)$. Here $\eta([0,s])$ is shown in dark gray, and $\eta([s,t])$ is colored light grey.  
    \textbf{Middle:} Another possible configuration.
    \textbf{Right:} Diagram for the definition of  {$(\widetilde \psi_s, \widetilde \eta_s)$} in Section~\ref{sec-radial-zipper}. Each of the quantum surfaces $\cC_j$  comes with four marked boundary points and a space-filling curve, as in Definition~\ref{def-cell}. We conformally weld $\cC_1$ (red) to $\cD_1 = (\bbD, \widetilde \psi_0, 0, 1)/{\sim_\gamma}$ along the two boundary arcs of $\cC_1$ adjacent to the starting point of its space-filling curve, to obtain $\cD_2$. Iterating this procedure (colors from red to purple, in order) gives $\cD_k$; we concatenate its curves and forget all marked points except the bulk point from $\cD_1$ (white) and the boundary endpoint of the curve of $\cC_k$ (purple), to get the quantum surface $(\bbD, \widetilde \psi_s, \widetilde \eta_s, 0, 1)/{\sim_\gamma}$. 
    Note that when each $\cC_j$ is conformally welded, {by construction} it will not ``wrap around'' the whole boundary of the other quantum surface. }
    \label{fig-lengths}
\end{figure}

\begin{corollary}\label{cor-area-law}
    For $\phi \sim \LF_\bbD^{(Q -\frac\gamma4, 0), (\frac{3\gamma}2, 1)}$ conditioned on having boundary length 1, the quantum area $\cA_\phi(\bbD)$ has the law of the inverse gamma distribution with shape parameter $\frac12$ and scale parameter $b = \frac18 \tan(\frac{\pi \gamma^2}8)$, i.e., the law of $\cA_\phi(\bbD)$ is 
    \[\mathds{1}_{a > 0} \sqrt{\frac {b}{\pi a^3} }e^{-\frac ba}  \, da.\]
\end{corollary}
\begin{proof}
The law of $X_t + Y_t$ is Brownian motion with quadratic variation $(2\mathbbm a \sin(\frac{\pi\gamma^2}8))^2\,dt = 4 \cot(\frac{\pi\gamma^2}8)\, dt$, and $\cA_\phi(\bbD)$ equals the hitting time of $-1$. The claim then follows from the well-known law of Brownian motion first passage times. 
\end{proof}

\begin{remark}\label{rem-lcft}
Corollary~\ref{cor-area-law}, together with the result \cite[Theorem 1.7]{RZ20b} and the computation of \cite[Section 4.4]{ARS21}, can be used to compute the correlation function of LCFT on the disk with a bulk insertion $\alpha = Q-\frac\gamma4$ and a boundary insertion $\beta = \frac{3\gamma}2$. This gives an alternative derivation of a special case of \cite[Theorem 1.2]{ARSZ23}, i.e., proves a special case of the physical proposal by \cite{hosomichi}.
\end{remark}

In Section~\ref{sec-radial-zipper} we define a \emph{radial quantum zipper} where, starting with a sample from $\LF_\bbD^{(Q + \frac\gamma4, 0), (-\frac\gamma2, 1)}$, we grow the quantum surface by conformal welding with independent quantum cells, giving rise to a coupling of LCFT with \emph{reverse} radial SLE. 
In Section~\ref{subsec-zipper} we prove Proposition~\ref{prop:radial-cell} in which we decorate $\LF_\bbD^{(Q + \frac\gamma4, 0), (\frac{3\gamma}2, 1)}$ by \emph{forward} radial SLE and look at the  quantum surfaces parametrized by the curve and its complement. Here, to switch between reverse and forward SLE,  we use the fact that for any fixed time, the curve generated by centered reverse radial SLE  has the law of forward radial SLE. 
In Section~\ref{sec-proof-radial}, since $\Delta_{Q +\frac\gamma4} = \Delta_{Q - \frac\gamma4}$ (with $\Delta_\alpha = \frac\alpha2 (Q - \frac\alpha2)$), we can use Girsanov's theorem  to obtain a variant of Proposition~\ref{prop:radial-cell} about $\LF_\bbD^{(Q - \frac\gamma4, 0), (\frac{3\gamma}2, 1)}$ (Proposition~\ref{prop:radial-finite}), and hence Theorem~\ref{thm-mot-disk-finite}.  

\subsection{A radial quantum zipper} \label{sec-radial-zipper}
Let $\kappa > 8$ and $\gamma = \frac4{\sqrt\kappa}$.  {The goal of this section is to prove  Lemma~\ref{lem:zipupdef}}, in which  we define and study a quantum zipper process {$(\psi_t, \eta_t)_{t \geq 0}$} where the marginal law of $\psi_0$ is $\LF_\bbD^{(Q + \frac\gamma4, 0), (-\frac\gamma2, 1)}$ and the time-evolution corresponds to conformally welding quantum cells to the boundary of the quantum surface.  {The proof of Lemma~\ref{lem:zipupdef} will depend on a result of \cite{Ang23} stated as Lemma~\ref{lem-chordal-zipper}.}

Let $ {\mathrm{BM}_\kappa}$ denote the law of (one-sided) correlated two-dimensional Brownian motion $(X_t, Y_t)_{t \geq 0}$ with $X_0 = Y_0=0$ and covariance given by~\eqref{eqn-thm:mot-chordal}.
Sample $(\widetilde \psi_0, (X_t, Y_t)) \sim \LF_\bbD^{(Q + \frac\gamma4, 0), (-\frac\gamma2, 1)} \times \mathrm{BM}_\kappa$, let $L_t = X_t + Y_t + \cL_{\widetilde \psi_0}(\partial \bbD)$, and let $\widetilde \tau$ be the first time $t$ that $L_t = 0$.  For $s \in (0,  \widetilde \tau)$ we define a random field and curve $(\widetilde \psi_s, \widetilde \eta_s)$ which correspond to ``zipping up for quantum time $s$'' as follows. See Figure~\ref{fig-lengths} (right). Choose finitely many times $0 =  s_1 < \dots < s_k = s$ such that for $j < k$ we have $(X_{s_{j}} - \inf_{u \in [s_{j}, s_{j+1}]} X_u)  + (Y_{s_{j}} - \inf_{u \in [s_{j}, s_{j+1}]} Y_u) < L_{s_{j}}$. For $j<k$ let $\cC_j = F((X_{\cdot + s_{j}} - X_{s_{j}}, Y_{\cdot + s_{j}} - Y_{s_{j}})_{[0, s_{j+1} - s_{j}]})$. We iteratively define quantum surfaces with the disk topology decorated by a bulk point, a boundary point, and a curve as follows. Let $\cD_1 = (\bbD, \widetilde \psi_0, 0, 1)/{\sim_\gamma}$, and iteratively for $j = 1, \dots, k-1$, we conformally weld $\cC_j$ to $\cD_{j}$ to obtain $\cD_{j+1}$. This is done by identifying the starting point of the curve of $\cC_j$ with the boundary point of $\cD_{j}$ and conformally welding the two boundary arcs of $\cC_j$ adjacent to this point to $\cD_{j}$ by quantum length {(this is possible since by assumption the two boundary arcs have total quantum length smaller than the quantum boundary length of $\partial \cD_j$)}.   Doing this $k-1$ times produces $\cD_k$, which we view as a quantum surface decorated by a bulk point (from {$\cD_1$}), a curve (obtained by concatenating the {$k-1$} curves from $\cC_1,\dots, {\cC_{k-1}}$), and a boundary point (the endpoint of the curve on the boundary). {We orient the curve so that it starts on the boundary and ends in the bulk of $\partial \cD_k$.} 
Finally, we conformally embed $\cD_k$ in $\bbD$, sending the bulk and boundary marked points to $0$ and $1$, to get $(\bbD, \widetilde \psi_s, \widetilde \eta_s, 0, 1)$. This gives our definition of $\widetilde \psi_s, \widetilde \eta_s$ for all $s < \widetilde \tau$; note that Lemma~\ref{lem-concat} implies this definition does not depend on the choice of $s_1, \dots, s_k$. 

{For each $s$, let $t(s)$ be the log conformal radius of $\bbD \backslash \widetilde \eta_s$ viewed from $0$, i.e., $t(s) = -\log |g'(0)|$ where $g:\bbD \to \bbD \backslash \eta_s$ is any conformal map fixing 0.} This gives a monotone reparametrization of the process which we denote by $(\psi_t, \eta_t)_{t \geq 0}$. We parameterize each curve $\eta_t:[0,t] \to \ol \bbD$ by log conformal radius, so $\eta_t(0) = 1$ and the  {conformal radius of $\bbD \backslash \eta_t([0,t'])$ viewed from 0 is $e^{-t'}$. {We first give a description of the process $(\psi_t, \eta_t)_{t \geq 0}$ in terms of the Liouville field and reverse SLE.} Recall $\bullet_\gamma$ from~\eqref{eq-coord-change}.}

\begin{lemma}\label{lem:zipupdef}
    For $\kappa>8$ and $\gamma=\frac{4}{\sqrt{\kappa}}$, let $M$ be the law of the process $(\psi_t,\eta_t)_{t\ge0}$ defined immediately above. Then 
    \begin{enumerate}[i)]
        \item For any a.s.\ finite stopping time $\tau$ for the filtration $\mathcal{F}_t$ generated by $(\eta_t)_{t\ge0}$, the law of $(\psi_\tau,\eta_\tau)$ is \\ $\LF_\bbD^{(Q+\frac{\gamma}{4},0),(-\frac{\gamma}{2},1)}\mathrm{rrSLE}_\kappa^\tau$, where $\mathrm{rrSLE}_\kappa^\tau$ denotes the law of centered reverse radial $\SLE_\kappa$ in $\bbD$ from 1 to 0 run until the stopping time $\tau$.
        \item For $0< t_1 < t_2$, let $\widetilde f_{t_1,t_2}:\bbD\to\bbD\backslash \eta_{t_2}([0,t_2-t_1])$ be the conformal map fixing 0 with $\widetilde f_{t_1,t_2}(1)=\eta_{t_2}(t_2-t_1)$, then $\psi_{t_1} = \widetilde{f}_{t_1,t_2}^{-1}\bullet_\gamma \psi_{t_2}$. 
    \end{enumerate}
\end{lemma}

 {
\begin{remark}
    \cite[Theorem 5.1]{ms-qle} constructed a process $(\psi_t, \eta_t)$ satisfying the conclusions of Lemma~\ref{lem:zipupdef}
    (but with the LQG field viewed modulo additive constant) by using a martingale argument to couple  GFF and SLE. For our purposes, however, we crucially require the mating-of-trees description of $M$ not present in \cite{ms-qle}. 

    The insertions $(\alpha, \beta) = (Q+\frac\gamma4, -\frac\gamma2)$ in the definition of $M$ satisfy $\alpha + \frac12\beta - Q = 0$, so the constant mode of the Liouville field has law $e^{(\alpha + \frac\beta2 - Q)c}dc = dc$ (up to multiplicative constant). The translation invariance of this law makes the Liouville field closely related to the GFF modulo additive constant, and hence the GFF/SLE coupling of \cite{ms-qle}. Moreover, the conformal invariance of the GFF modulo additive constant is the underlying reason why prefactors cancel in our subsequent argument (below~\eqref{eq:chordal-zipup}).
\end{remark}
}


\begin{proof}[Proof of Lemma~\ref{lem:zipupdef}]
From the definition of $M$, 
$(\bbD, \psi_{t_2}, 0, 1)/{\sim_\gamma}$ is obtained from conformally welding $(\bbD, \psi_{t_1}, 0, 1)/{\sim_\gamma}$ with another quantum surface, so $(\bbD \backslash \eta_{t_2}([0, t_2 - t_1]), \psi_{t_2}, 0, \eta_{t_2}(t_2 - t_1))/{\sim_\gamma} = (\bbD, \psi_{t_1}, 0, 1)/{\sim_\gamma}$. This gives $\psi_{t_1} = \widetilde f_{t_1, t_2}^{-1} \bullet_\gamma \psi_{t_2}$ so ii) holds. 

For i), we first {apply a change of coordinates from $(\bbD, 1, -1)$ to $(\bbH, 0, \infty)$ to change the radial process $(\psi_t, \eta_t)_{t \geq 0}$ into a chordal process $(\hat \psi_t, \hat \eta_t)_{t \geq 0}$ in $(\bbH, 0, \infty)$,} apply Lemma~\ref{lem-chordal-zipper} for the chordal process in $\bbH$, and finally convert back to the radial process in $\bbD$.


For a sample $(\psi_t, \eta_t)_{t \geq 0} \sim M$, let $\tau_0$ be the time $t$ that $\widetilde f_{0,t}(-1) = 1$, or in other words the time the boundary point $p_0 = -1$ of $(\bbD, \phi_0)$ intersects the zipped-in region (colored region in Figure~\ref{fig-lengths} (right)). Let $g_0: \bbD \to \bbH$ be the conformal map such that $g_0(0)= i$ and $g_0(1) = 0$. For $t < \tau_0$ let $p_t = \widetilde f_{0,t}(p_0) \in \partial \bbD \backslash \{ 1\}$, and let $g_t : \bbD \to \bbH$ be the conformal map such that $g_t(1) = 0, g_t(p_t) = \infty$, and $(g_t \circ \widetilde f_{0, t} \circ g_0^{-1})(z) = z + O(1)$ as $z \to \infty$. This gives us a process $(g_t \bullet_\gamma \psi_t, g_t \circ \eta_t)_{[0, \tau_0)}$ of (field, curve) pairs in $\bbH$; we reparametrize time to obtain a process $(\hat \psi_t, \hat \eta_t)_{[0, \infty)}$ such that the half-plane capacity of the trace of $\hat \eta_t$ is $2t$, and $\hat \eta_t:[0,t] \to \ol \bbH$ is parametrized by half-plane capacity. 
By Lemma~\ref{lem:lfcoord}, the law of $(\hat \psi_0, (X_t, Y_t)_{t\ge0})$ is $\LF_\bbH^{(Q + \frac\gamma4, i), (-\frac\gamma2, 0)} \times \mathrm{BM}_\kappa$, and by our choice of $g_t$ the conformal maps $\widetilde f_{\bbH,t}: \bbH \to \bbH \backslash \hat \eta_t([0,t])$ satisfying $\widetilde f_{\bbH,t} (0) = \hat \eta_t(t)$ and $\widetilde f_{\bbH,t}(z) = z + O(1)$ as $z \to \infty$ also satisfy $\hat \psi_0 = \widetilde f_{\bbH,t}^{-1} \bullet_\gamma \hat \psi_t$. 

 {
By Lemma~\ref{lem-concat} the process $(\hat \psi_t, \hat \eta_t)$ is as described in Section~\ref{sec-quantum-zipper}, so by Lemma~\ref{lem-chordal-zipper},}
for any stopping time $\sigma$ for the filtration $\cF_t = \sigma(\hat\eta_t)$, the law of $(\hat\psi_{\sigma},\hat\eta_{\sigma})$ is
\begin{equation}\label{eq:chordal-zipup}
\frac{1}{\mathcal{Z}_{Q+\frac\gamma4}(\widetilde{f}_{\bbH, \sigma}(i))}\LF_\bbH^{(Q+\frac{\gamma}{4},\widetilde f_{\bbH, \sigma}(i)), (-\frac{\gamma}{2},0)}\,\mathrm{rSLE}_{\kappa,\kappa+6}^{\sigma}, \qquad \mathcal{Z}_{Q+\frac\gamma4}(z) = (2\mathrm{Im}\,z)^{-\frac{(\kappa+6)^2}{8\kappa}} |z|^{\frac{\kappa+6}{{\kappa}}}.
\end{equation} 
Applying the conformal map $\bbH \to \bbD$ sending $\widetilde f_\sigma (i)$ to $0$ and $0$ to $1$, and using the LCFT change of coordinates Lemma~\ref{lem:lfcoord} and reverse SLE change of coordinates Lemma~\ref{lem:rslecoord}, we obtain i) for any stopping time $\tau$ with $\tau \leq \tau_0$ a.s.. (Note that the prefactor incurred from Lemma~\ref{lem:lfcoord} cancels with the factor $1/\mathcal Z_{Q+\frac\gamma4}(\widetilde f_{\bbH, \sigma}(i))$ from~\eqref{eq:chordal-zipup}).

As we will see,  the above result can be iterated to get i) for all $\tau$. By the previous paragraph, the law of $(\psi_{\tau_0}, \eta_{\tau_0})$ is $\LF_\bbD^{(Q + \frac\gamma4, 0), (-\frac\gamma2, 1)} \times \mathrm{rrSLE}_\kappa^{\tau_0}$, so by the Markov property of Brownian motion, the law of $((\psi_{\tau_0 + t}, \eta_{\tau_0 + t}|_{[0, t]})_{t \geq 0}, \eta_{\tau_0})$ is $M \times \mathrm{rrSLE}_\kappa^{\tau_0}$. Define $\tau_1$ for $(\psi_{\tau_0 + t}, \eta_{\tau_0 + t}|_{[0, t]})_{t \geq 0}$ in the same way that $\tau_0$ was defined for $(\psi_t, \eta_t)_{t \geq 0}$, so $\tau_0, \tau_1$ are i.i.d.. 
Conditioning on $\eta_{\tau_0}$ and applying the result of the previous paragraph, we see that i) holds for any stopping time $\tau \leq \tau_0 + \tau_1$. Proceeding iteratively,  we may define $\tau_k$ for all $k$, and i) holds for all $\tau \leq \sum_{i \leq k} \tau_i$. Since  the $\tau_k$ are i.i.d.\ positive random variables we have $\sum_k \tau_k \to \infty$ a.s., completing the proof of i). 
\end{proof}

The following lemma essentially tells us that if we run the process $(\psi_t, \eta_t)_{t \geq 0}$ until a random amount of quantum area has been added, if the added region is simply connected then it parametrizes a quantum cell independent of $\psi_0$.

The following lemma is the radial analog of \cite[Proposition 5.7]{Ang23}. 
\begin{lemma}\label{prop:disk-zip-up}
   Let $\kappa>8$ and $\gamma=\frac{4}{\sqrt\kappa}$. Sample $((\psi_t,\eta_t)_{t\ge0},A)$ from $M\times \mathds{1}_{a>0}da$. Restrict to the event that there is a time $\tau>0$ such that $\mathcal{A}_{\psi_\tau}(\eta_\tau([0,\tau]))=A$. 
   Then the law of $(\psi_\tau,\eta_\tau,\tau)$ is $C\cdot\LF_\bbD^{(Q+\frac{\gamma}{4},0),(\frac{3\gamma}{2},1)}\times\mathrm{raSLE}_\kappa^t\mathds{1}_{t>0}dt$ for some constant $C>0$. Here $\mathrm{raSLE}_\kappa^t$ denotes the law of radial $\SLE_\kappa$ in $\bbD$ from 1 to 0 parametrized by log-conformal radius stopped at time $t$.
\end{lemma}
\begin{proof}
Here is a proof sketch;  {just for this proof, we use the shorthand $s_+ := \max\{s,0\}$}. First, if we fix $\delta > 0$ and sample $((\psi_t, \eta_t)_{t \geq 0}, A, T)\sim \delta^{-1} \mathds{1}_{T \in [\tau, \tau + \delta]} M \times \mathds{1}_{A>0} \, dA \times dT$ then the marginal law of $((\psi_t, \eta_t)_{t \geq 0}, A)$ is $M\times \mathds{1}_{a>0}da$, so the marginal law of $(\psi_\tau, \eta_\tau, \tau)$ is the same as in Lemma~\ref{prop:disk-zip-up}.  {In this new setup, the constraint $\{T \in [\tau, \tau + \delta]\} = \{\tau\in[(T-\delta)_+,T]\}$ is the same as $$\cA_{\psi_T}(\eta_T([0,T]))\geq A\geq \cA_{\psi_{(T-\delta)_+}}(\eta_{(T-\delta)_+}([0,(T-\delta)_+])).$$ Note that the lower bound equals $\cA_{\psi_T}(\eta_T([\delta\wedge T,  T]))$ by ii) of Lemma~\ref{lem:zipupdef}, so using i) of Lemma~\ref{lem:zipupdef}, the law of $(A,\psi_T,\eta_T,T)$ is then
$$\delta^{-1}\mathds{1}_{a \in [\cA_\psi(\eta([\delta\wedge t,  t])),\cA_\psi(\eta([0,t]))]}da\times \LF_\bbD^{(Q + \frac\gamma4, 0), (-\frac\gamma2, 1)}(d\psi)\times \mathrm{rrSLE}_\kappa^t(d\eta) \mathds{1}_{t > 0} \, dt.$$}
Let $z = \eta_T(T - \tau)$. Since  {$\mathcal{A}_{\psi_T}(\eta_T([T-\tau,T])) = A$, we have} $\{T \in [\tau, \tau + \delta]\} = \{ z \in \eta_T([0,\delta \wedge T])\}$.   {On the other hand, from the definition of $\tau$, $z$ can be viewed as a point sampled on $\eta_T([0,\delta \wedge T])$ according to the measure $\cA_{\psi_T}$. Therefore} the law of $(z, \psi_T, \eta_T, T)$ is 
$ {\delta^{-1}}\mathds{1}_{z \in \eta([0,\delta\wedge t])}\cA_{\psi}(dz)\LF_\bbD^{(Q + \frac\gamma4, 0), (-\frac\gamma2, 1)}(d\psi)\times \mathrm{raSLE}_\kappa^t(d\eta) \mathds{1}_{t > 0} \, dt$. 
Note we have obtained the term $\mathrm{raSLE}_\kappa^t \mathds{1}_{t>0}\,dt$ using the symmetry between forward and reverse radial $\SLE_\kappa$ at fixed time $t$.
Using Lemma~\ref{lem:qtypical}, this law is \[ {\delta^{-1}}\LF_\bbD^{(Q + \frac\gamma4, 0), (-\frac\gamma2, 1), (\gamma, z)}(d\psi) \mathds{1}_{z \in \eta([0,\delta\wedge t])} dz \, \mathrm{raSLE}_\kappa^t(d\eta) \mathds{1}_{t > 0} \, dt.\]
As $\delta \to 0$ we have $T - \tau \to 0$ so $z \to 1$, so in the limit the field has the singularity $\gamma G(\cdot, 1) - \frac\gamma4 G(\cdot, 1) = \frac12 (\frac{3\gamma}2) G(\cdot, 1)$ at $1$. This explains the term $\LF_\bbD^{(Q + \frac\gamma4, 0), (\frac{3\gamma}2, 1)}$. The main difficulty in this argument is in taking limits of infinite measures; this is done by truncating on finite events and taking limits of finite measures. 

The argument outlined above is implemented in the proof of 
\cite[Proposition 5.7]{Ang23}, a chordal analog of our desired result; we refer the reader there for details.  The only part of that proof that does not immediately carry over to our setting is a certain finiteness claim \cite[Lemma 5.8]{Ang23}, whose analog in our setting can be stated as follows. 
For $\rho$ the uniform probability measure on $\{ z: |z| = \frac12\}$ (the precise choice of $\rho$ is unimportant), we have 
    \eqb\label{eq-EN}
(M\times \mathds{1}_{a>0} da)[E_N]< \infty \text{ where }   E_N :=  \{\tau, |(\psi_0, \rho)|, |(\psi_\tau, \rho)| < N \}.  
    \eqe
    Given this, the proof of our Lemma~\ref{prop:disk-zip-up} is identical to that of \cite[Proposition 5.7]{Ang23}. Thus it suffices to prove~\eqref{eq-EN}.

    First, we observe $(M \times \mathds{1}_{a>0}da)[E_N] \leq (M \times \mathds{1}_{a>0}da)[\widetilde E_N] = M[\cA_{\psi_N}(\eta_N([0,N])) \mathds{1}_{|(\psi_0, \rho)| < N}]$ where $\widetilde E_N = \{\tau , |(\psi_0, \rho)|<N\}$. Now, {our choice of parametrization implies the conformal radius of $\eta_N([0,N])$ viewed from $0$ is $e^{-N}$}, so the Koebe quarter theorem implies that the ball $B_{e^{-N}/4}(0)$ is contained in $\bbD \backslash \eta_N([0,N])$. By Lemma~\ref{lem:zipupdef} the $M$-law of $(\psi_N, \eta_N)$ is $\LF_\bbD^{(Q + \frac\gamma4, 0), (-\frac\gamma2, 1)} \times \mathrm{raSLE}_\kappa^N$, so it suffices to show the finiteness of 
    \eqb \label{eq-EN-2}
(\LF_\bbD^{(Q + \frac\gamma4, 0), (-\frac\gamma2, 1)} \times \mathrm{raSLE}_\kappa^N)[\cA_\psi(\bbD \backslash B_{e^{-N}/4}(0)) \mathds{1}_{|({ f_N} \bullet_\gamma \psi, \rho)|< N}],
    \eqe
    where for  the $\mathrm{raSLE}_\kappa^N$ curve $\eta$ the conformal map ${ f_N}: \bbD\backslash \eta([0,N]) \to \bbD$ satsifies $ f_N(0) =0$ and $ f_N(\eta(N)) = 1$. Writing $\bbE$ to denote expectation with respect to $(h, \eta) \sim P_\bbD\times \mathrm{raSLE}_\kappa^N$ and $\widetilde h = h + (Q + \frac\gamma4) G_\bbD(\cdot, 0) -\frac\gamma4 G_\bbD(\cdot, 1)$, this equals 
    \begin{align*}
        \bbE\big[\int_\bbR \cA_{\widetilde h + c}(\bbD \backslash B_{e^{-N}/4}(0)) \mathds{1}_{|( f_N \bullet_\gamma \widetilde h, \rho) + c|<N}\,dc\big] &= \bbE \big[\int_{-( f_N \bullet_\gamma \widetilde h, \rho)-N}^{-( f_N \bullet_\gamma \widetilde h, \rho)+N} e^{\gamma c} \cA_{\widetilde h}(\bbD \backslash B_{e^{-N}/4}(0)) \,dc \big] \\
        &=\frac1\gamma(e^{\gamma N} - e^{-\gamma N}) \bbE \big[e^{-\gamma ( f_N \bullet_\gamma \widetilde h, \rho)} \cA_{\widetilde h}(\bbD \backslash B_{e^{-N}/4}(0)) \big].
    \end{align*}
    To see this is finite, first note that $Z := \bbE[e^{-\gamma( f_N \bullet_\gamma \widetilde h, \rho)}]<\infty$ by standard conformal distortion estimates. 
    Next, by Girsanov's theorem, the expression equals $\frac1\gamma (e^{\gamma N} - e^{-\gamma N})Z\bbE[\cA_{\hat h} (\bbD \backslash B_{e^{-N}/4}(0)]$ where $\hat h = h + (Q + \frac\gamma4) G_\bbD(\cdot, 0) -\frac\gamma4 G_\bbD(\cdot, 1) - \gamma \int G_\bbD(\cdot, w) (( f_N^{-1})_* \rho)(dw)$. To finish, we note that $\hat h - h$ is bounded above by a constant on $\bbD \backslash B_{e^{-N}/4}(0)$, and that $\bbE[\cA_h(\bbD \backslash B_{e^{-N}/4}(0))]<\infty$ by standard GMC moment results, see for instance \cite[Proposition 3.5]{RV10}. We conclude that~\eqref{eq-EN-2}, and hence~\eqref{eq-EN}, is finite.
\end{proof}

Finally, between two ``quantum typical'' times for $(\psi_t, \eta_t) \sim M$, given the field and curve at the earlier time, on the event the zipped-in quantum surface is simply connected, it is a quantum cell with  a boundary length restriction. 
\begin{lemma}\label{lem:disk-zip-up}
         Let $\kappa>8$ and $\gamma=\frac{4}{\sqrt\kappa}$, and fix $a_1, a_2 > 0$. Sample $(\psi_t,\eta_t)_{t\ge0}$ from $M$ and restrict to the event that there is a time $\tau_2>0$ such that $\mathcal{A}_{\psi_{\tau_2}}(\eta_{\tau_2}([0,\tau_2]))=a_1 +a_2$. Let $\tau_1$ be the time that $\mathcal A_{\psi_{\tau_1}}(\eta_{\tau_1}([0,\tau_1])) = a_1$. Conditioned on $(\psi_{\tau_1}, \eta_{\tau_1})$, the law of $(\eta_{\tau_2}([0,\tau_2 - \tau_1]), \psi_{\tau_2},\eta_{\tau_2} |_{[0,\tau_2 - \tau_1]})/{\sim_\gamma}$ restricted to the event $\{\eta_{\tau_2}([0,\tau_2 - \tau_1]) \textup{ is simply connected}\}$ is 
    \[
        \mathds{1}_{X_{a_2}^+(\cC) + Y_{a_2}^+(\cC) < \cL_{\psi_{\tau_1}}(\partial \bbD)} P_{a_2}(d \cC)
    \]
    where $X_{a_2}^+$ and $Y_{a_2}^+$ are as in Definition~\ref{def-cell}.
     \end{lemma}
\begin{proof}
{
Let $(X_t, Y_t)_{t \geq 0}$ be the process in the definition of $M$, then the law of  $\widetilde \cC := F((X_{\cdot+a_1}, Y_{\cdot + a_1}))_{[0,a_2]}$ is $P_{a_2}$, and reversing the orientation of the curve of $\widetilde \cC$  gives $\cC := (\eta_{\tau_2}([0,\tau_2 - \tau_1]), \psi_{\tau_2},\eta_{\tau_2} |_{[0,\tau_2 - \tau_1]})/{\sim_\gamma}$. 
By construction $\{ \eta_{\tau_2}([0,\tau_2 - \tau_1]) \text{ is simply connected}\} = \{X_{a_2}^-(\widetilde \cC) + Y_{a_2}^-(\widetilde \cC) < \cL_{\psi_{\tau_1}}(\partial \bbD) \}$, and since  $X_{a_2}^- (\widetilde \cC) = Y_{a_2}^+(\cC)$ and $Y_{a_2}^-(\widetilde \cC) = X_{a_2}^+(\cC)$, this event equals $ \{X_{a_2}^+( \cC) + Y_{a_2}^+( \cC) < \cL_{\psi_{\tau_1}}(\partial \bbD) \} $ as needed.}
\end{proof}

\subsection{Cutting an infinite volume LCFT disk until a quantum typical time}\label{subsec-zipper}
The aim of this section is to prove Proposition~\ref{prop:radial-cell} below. We write $\mathrm{raSLE}_\kappa^t$ for the law of radial $\SLE_\kappa$ in $\bbD$ from 1 to 0 stopped at time $t$, {and $\mathrm{raSLE}_\kappa^z$ for the law of radial $\SLE_\kappa$ in $\bbD$ from 1 to 0 stopped when it hits $z \in \ol \bbD \backslash \{0\}$}.


\begin{proposition}\label{prop:radial-cell}
Suppose $\kappa > 8$ and $\gamma = \frac4{\sqrt\kappa}$.
    Sample $(\phi,\eta,A)$
    from the measure 
    \eqb\label{eq-no-t'}
    \LF_\bbD^{(Q+\frac{\gamma}{4},0),(\frac{3\gamma}{2},1)}\times \mathrm{raSLE}_\kappa\times \mathds{1}_{a>0}da
    \eqe and parametrize  {$\eta$} by its  {$\mathcal{A}_{\phi}$} quantum area. For $a\ge0$,  let $f_a:\bbD\backslash \eta([0,a])\to\bbD$ be the conformal map such that $f_a(0) = 0$ and $f_a(\eta(a))=1$. Let $\phi_a = f_a\bullet_\gamma\phi$,  {$\widetilde\eta_a=f_a\circ\eta|_{[a,\infty)}$,} 
     and $\mathcal{C}_a =  (\eta([0,a]), \phi, \eta|_{[0,a]})/{\sim_\gamma}$. Then the law of $(\phi_A, \widetilde\eta_A,A)$ is given by\footnote{ {There is a slight abuse of notation here: the curve $\widetilde \eta_A$ should be viewed as parametrized by log-conformal radius rather than by quantum area for~\eqref{eq:radial-cell0} to hold. We do this because this section is already notationally dense.}}  
    \eqb\label{eq:radial-cell0}
    \LF_\bbD^{(Q+\frac{\gamma}{4},0),(\frac{3\gamma}{2},1)} \times \mathrm{raSLE}_\kappa \times  \mathds{1}_{a>0}da.
    \eqe
    Moreover, the law of $(\phi_A, \widetilde \eta_A, \mathcal{C}_A,A)$ restricted to the event that $\eta([0,A])$ is simply connected is given by 
    \begin{equation}\label{eq:radial-cell}
        \mathds{1}_{X_a^+(\cC_a) + Y_a^+(\cC_a) < \cL_{\phi_a}(\partial \bbD)}\LF_\bbD^{(Q+\frac{\gamma}{4},0),(\frac{3\gamma}{2},1)}(d\phi_a)\times \mathrm{raSLE}_\kappa\times P_a(d\cC_a)\, \mathds{1}_{a>0}da.
    \end{equation}
    where $X_a^+, Y_a^+$ are as in Definition~\ref{def-cell}. 
\end{proposition}

\begin{figure}
    \centering
    \includegraphics[scale=0.8]{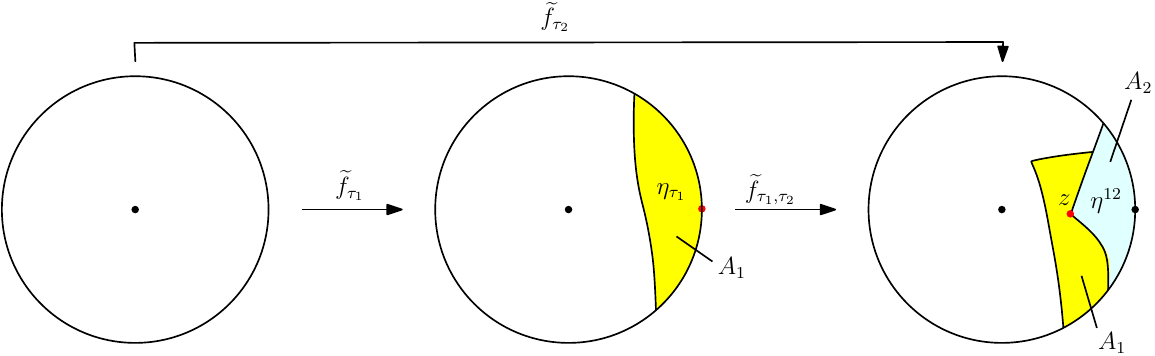}
    \caption{ {Setup for the proof of Proposition~\ref{prop:radial-cell}. We sample $A_1,A_2$ from $\mathds{1}_{A_1,A_2>0}dA_1dA_2$. The middle panel corresponds to the time $\tau_1$ where the yellow quantum cell (filled with the curve $\eta_{\tau_1}$) has quantum area $A_1$ has been ``zipped in'', while in the right panel we continue to time $\tau_2$ where we have conformally welded the blue quantum cell (filled with the curve $\eta^{12}$) with quantum area $A_2$. In the right cell, the curve $\eta_{\tau_2}$ is the concatenation of the curves in the blue and yellow cells, 
    and $z$ corresponds to the point $\eta_{\tau_2}(\tau_2-\tau_1)$.}}
    \label{fig:pf-prop:radial-cell}
\end{figure}

\begin{proof}[Proof of Proposition~\ref{prop:radial-cell}]
 {Here is a proof outline. Consider the setup in Figure~\ref{fig:pf-prop:radial-cell} where we conformally weld quantum cells of areas $A_1$ and $A_2$ in the definition of $M$. Since $A_1$, $A_2$ are sampled from Lebesgue measure, the point $z$ can be viewed as being sampled from $\cA_{\psi_{\tau_2}}$-measure on $\eta_{\tau_2}([0,\tau_2])$, and the radial $\SLE_\kappa$ $\eta_{\tau_2}$ can be decomposed into $\wt{f}_{\tau_1,\tau_2}\circ\eta_{\tau_1}$ and $\eta^{12}$. Then by Lemma~\ref{prop:disk-zip-up}   $(\psi_{\tau_2}, \eta^{12}, A_2)\stackrel d=  (\phi,\eta|_{[0,A]},A)$. Passing to the middle panel via $\wt{f}_{\tau_1,\tau_2}^{-1}$ gives a description of the law of $(\phi_A,\eta_A,A)$ in terms of that of $(\psi_{\tau_1},\eta_{\tau_1},\tau_1)$, and the conclusion follows from another application of Lemma~\ref{prop:disk-zip-up}.}

To streamline notation in this proof, we will often use the same notation for a random object as in the description of its law (in the indented equations), or similar notation (e.g.\ use $d\psi_{t_2}$ in a description of the law of $\psi_{\tau_2}$).
    To begin with, sample $(\{(\psi_t,\eta_t)_{t\ge0}\},A_1,A_2)$ from $M\times \mathds{1}_{A_1,A_2>0}dA_1dA_2$, and let $\tau_1$ (resp.\ $\tau_2$) be the time $t$ when $\mathcal{A}_{\psi_t}(\eta_t([0,t]))$ equals $A_1$ (resp.\ $A_1+A_2$). {We restrict to the event $E$ that these times exist ($\tau_1 < \tau_2 < \infty$).} Let $z = \eta_{\tau_2}(\tau_2-\tau_1)$, $S=A_1+A_2$, and $\eta^{12} = \eta_{\tau_2}|_{[0,\tau_2-\tau_1]}$. Then the law of $((\psi_t,\eta_t)_{t\ge0},A_1,S)$ is ${\mathds{1}_E}M\times \mathds{1}_{A_1\in[0,S]}dA_1\mathds{1}_{S>0}dS$, so by Lemma~\ref{prop:disk-zip-up}  {applied to $((\psi_t,\eta_t)_{t\ge0},S)$}, the law of $(A_1,\psi_{\tau_2},\eta_{\tau_2},\tau_2)$ is
    $$ C\cdot \mathds{1}_{A_1\in[0,\mathcal{A}_{\psi_{t_2}}(\eta_{t_2}([0,t_2])]}dA_1\,\LF_\bbD^{(Q+\frac{\gamma}{4},0),(\frac{3\gamma}{2},1)}(d\psi_{t_2})\,\mathrm{raSLE}_\kappa^{t_2}(d\eta_{t_2})\,\mathds{1}_{t_2>0}dt_2.$$
    Since $z$ is the point where $\eta_{\tau_2}$ covers $S-A_1$ units of quantum area when hitting $z$, it follows that the law of $(z,\psi_{\tau_2}, \eta_{\tau_2},\tau_2)$ is
    $$C\cdot \mathds{1}_{z\in \eta_{t_2}([0,t_2])}\mathcal{A}_{\psi_{t_2}}(dz)\, \LF_\bbD^{(Q+\frac{\gamma}{4},0),(\frac{3\gamma}{2},1)}(d\psi_{t_2})\,\mathrm{raSLE}_\kappa^{t_2}(d\eta_{t_2})\,\mathds{1}_{t_2>0}dt_2.$$
    Then by   Lemma~\ref{lem:sledmp} below, the law of $(z,\psi_{\tau_2}, \eta^{12},(\eta_{\tau_1}, \tau_1))$ is 
    \eqb\label{eq-5-tuple}
    \mathcal{A}_{\psi_{t_2}}(dz)\LF_\bbD^{(Q+\frac{\gamma}{4},0),(\frac{3\gamma}{2},1)}(d\psi_{t_2})\, \mathrm{raSLE}_\kappa^{z}(d\eta^{12})\times[ C\cdot \mathrm{raSLE}_\kappa^{t_1}(d\eta_{t_1})\,\mathds{1}_{t_1>0} dt_1]
    \eqe
    {where $\mathrm{raSLE}_\kappa^z$ is as defined before Proposition~\ref{prop:radial-cell}.}

 Since $(\phi,\eta,A)$ is sampled from~\eqref{eq-no-t'} and $\eta$ is parametrized by quantum area, \ {the law of $(\eta(A),\phi,\eta|_{[0,A]})$ is 
 $\cA_{\phi}(du)\LF_\bbD^{(Q+\frac{\gamma}{4},0),(\frac{3\gamma}{2},1)}(d\phi)\, \mathrm{raSLE}_\kappa^{u}(d\eta)$}. 
 Then, by the domain Markov property of radial $\SLE_\kappa$, if we  {instead sample $(\phi, \eta, A, t')$ from 
 \eqb\label{eq-indep-t'}
 \LF_\bbD^{(Q+\frac{\gamma}{4},0),(\frac{3\gamma}{2},1)}(d\phi)\times \mathrm{raSLE}_\kappa(d\eta)\times \mathds{1}_{a>0}\, da \times  [C \mathds{1}_{t > 0} \,dt]
 \eqe
 (or ``independently sample $t'$ from $[C\mathds{1}_{t>0} dt]$'')} then the law of $(\eta(A),\phi,\eta|_{[0,A]}, ({\widetilde\eta_A|_{[0,t']}},t'))$ 
 is
 {
\[ \cA_{\phi}(du)\LF_\bbD^{(Q+\frac{\gamma}{4},0),(\frac{3\gamma}{2},1)}(d\phi)\, \mathrm{raSLE}_\kappa^{u}(d\eta) \times [C \cdot \mathrm{raSLE}_\kappa^{t} \mathds{1}_{t > 0} \, dt].\]
 This law agrees with~\eqref{eq-5-tuple} up to renaming random variables}, so $(\eta(A),\phi,\eta|_{[0,A]}, {\widetilde\eta_A|_{[0,t']}},t') \stackrel d= (z,\psi_{\tau_2}, \eta^{12}, \eta_{\tau_1},\tau_1)$.
 Since $A_2=\mathcal{A}_{\psi_{\tau_2}}(\eta^{12})$, $\psi_{\tau_1} = \widetilde{f}_{\tau_1,\tau_2}^{-1}\bullet_\gamma\psi_{\tau_2}$ where $\widetilde{f}_{\tau_1,\tau_2}:\bbD\to\bbD\backslash\eta^{12}$ is the conformal map fixing 0 and sending 1 to the tip of $\eta^{12}$, it follows that $(\phi_A,A,\widetilde\eta_A|_{[0,{t'}]},t') \stackrel d= (\psi_{\tau_1},A_2, \eta_{\tau_1}, \tau_1)$. 

    On the other hand, by Lemma~\ref{prop:disk-zip-up}, the law of {$(\psi_{\tau_1},A_2, \eta_{\tau_1},\tau_1)$} is 
    \eqb\label{eq-time-a1}
    \LF_\bbD^{(Q+\frac{\gamma}{4},0),(\frac{3\gamma}{2},1)}(d\psi_{t_1})\,\mathds{1}_{A_2>0}dA_2\times \mathrm{raSLE}_\kappa^{t_1}(d\eta_{t_1})\, [C \mathds{1}_{t_1>0} dt_1].
    \eqe
    Note the term $[C \mathds{1}_{t_1>0} dt_1]$ above corresponds to $[C\mathds{1}_{t>0}dt]$ in~\eqref{eq-indep-t'}, so by varying $t'$, for  $(\phi, \eta, A)$ sampled from~\eqref{eq-no-t'} the law of $(\phi_A, \widetilde\eta_A,A)$  is given by~\eqref{eq:radial-cell}. 
    This concludes the proof of the first claim. 

    {
    For the second claim, we repeat the above except we restrict to the event $F := \{ \eta^{12} \text{ is simply connected}\}$ throughout. Then the law of $(z,\psi_{\tau_2}, \eta^{12},\eta_{\tau_1}, \tau_1)$ is $\mathds{1}_F$ times~\eqref{eq-5-tuple}, and by Lemmas~\ref{prop:disk-zip-up} and~\ref{lem:disk-zip-up}, the law of $(\psi_{\tau_1}, \eta_{\tau_1}, (\eta^{12}([0,\tau_2 - \tau_1]), \psi_{\tau_2}, \eta^{12})/{\sim_\gamma}, A_2, \tau_1)$ is 
    \[\mathds{1}_{X_{A_2}^+(\cC) + Y_{A_2}^+(\cC) < \cL_{\psi_{t_1}}(\partial \bbD)} \LF_\bbD^{(Q+\frac{\gamma}{4},0),(\frac{3\gamma}{2},1)}(d\psi_{t_1})\,\mathrm{raSLE}_\kappa^{t_1}(d\eta_{t_1})\,P_{A_2}(d\cC) \, \mathds{1}_{A_2>0}dA_2\, C \mathds{1}_{t_1>0} dt_1, \]
    c.f.~\eqref{eq-time-a1}. The same argument as that of the first claim then gives the second claim. }
\end{proof}

\begin{figure}
    \centering
    \includegraphics[scale=0.8]{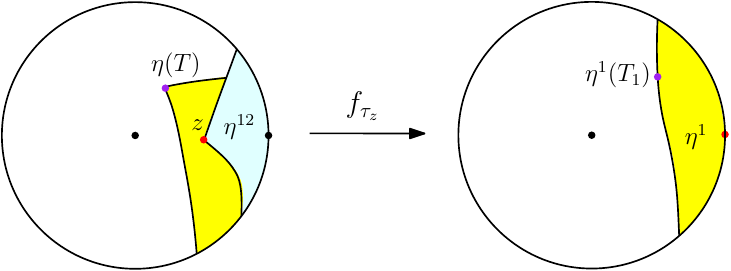}
    \caption{An illustration of Lemma~\ref{lem:sledmp}, where $z$ is fixed, $T$ is sampled from Lebesgue measure restricted to the event that the radial $\SLE_\kappa$ curve $\eta$ covers $z$ at time  $\tau_z<T$. Let $T_1=T-\tau_z$ and $f_{\tau_z}$ be the centered Loewner map at time $\tau_1$. We show that the law of {$(\eta^{12},(\eta^1,T_1))$ is $\mathrm{raSLE}_\kappa^z\times [\mathrm{raSLE}_\kappa^{t_1}\,\mathds{1}_{t_1>0}dt_1]$}.
    }
    \label{fig:sledmp}
\end{figure}

In the above proof we needed the following lemma,  {see Figure~\ref{fig:sledmp}.}

\begin{lemma}\label{lem:sledmp}
    Fix $z\in\bbD$, and sample $(\eta,T)$  from $\mathds{1}_{z\in\eta([0,t])}\mathrm{raSLE}_\kappa(d\eta)\mathds{1}_{t>0}\,dt$  {where $\eta$ is parametrized by log-conformal radius seen from 0.} Let $\tau_z$ be the time when $\eta$ hits $z$, $T_1=T-\tau_z$ and $\eta^{12}=\eta|_{[0,\tau_z]}$. Let  {$f_{\tau_z}:\bbD\backslash\eta([0,\tau_z])\to\bbD$} be the centered Loewner map of $\eta$ at time $\tau_z$, and  {$\eta^1 = f_{\tau_z}\circ\eta(\cdot+\tau_z)|_{[0,T_1]}$}. 
    Then the law of {$(\eta^{12},(\eta^1,T_1))$ is $\mathrm{raSLE}_\kappa^z\times [\mathrm{raSLE}_\kappa^{t_1}\,\mathds{1}_{t_1>0}dt_1]$}, where  $\mathrm{raSLE}_\kappa^z$ is the law of radial $\SLE_\kappa$ run until it hits $z$,  {and $\mathrm{raSLE}_\kappa^{t_1}$ is the law of the radial $\SLE_\kappa$ curve stopped at the time when the log-conformal radius seen from 0 equals $t_1$ as in Lemma~\ref{prop:disk-zip-up}}.
\end{lemma}
\begin{proof}
By a change of variables, the law of $(\eta^{12}, T_1)$ is $\mathrm{raSLE}_\kappa^z(d\eta^{12})\times \mathds{1}_{t_1>0}dt$. By the domain Markov property of radial $\SLE$, conditioned on $\eta^{12}$ and $T_1$, the law of $\eta^1$ is $\mathrm{raSLE}_\kappa^{T_1}$. This finishes the proof.
\end{proof}

\subsection{Proof of Theorem~\ref{thm-mot-disk-finite}}\label{sec-proof-radial}
In this section we prove Theorem~\ref{thm-mot-disk-finite}. We first use Proposition~\ref{prop:radial-cell} about $\LF_\bbD^{(Q + \frac\gamma4, 0), (\frac{3\gamma}2, 1)}$ to obtain an analogous result for $\LF_\bbD^{(Q - \frac\gamma4, 0), (\frac{3\gamma}2, 1)}$ (Proposition~\ref{prop:radial-finite}). The idea is to weight the field to change $Q + \frac\gamma4$ into $Q - \frac\gamma4$ via the Girsanov theorem; this is done in Lemmas~\ref{lem-weight} and~\ref{lem-weight-fancy}.

Let $B_\eps(0) := \{ z \in \bbC \:: \: |z|< \eps\}$, and let $\theta_\eps$ denote the uniform probability measure on the circle $\partial B_\eps(0)$. 
\begin{lemma}\label{lem-weight}
Let $\alpha_1, \alpha_2, \beta \in \bbR$ and $\eps \in (0,1)$. Let $\widetilde\LF_{\bbD, \eps}^{(\alpha_2, 0), (\beta, 1)}$ be the law of $\psi|_{\bbD \backslash B_\eps(0)}$, where $\psi \sim \LF_{\bbD}^{(\alpha_2, 0), (\beta, 1)}$. 
     Sample $\phi$ from the measure $\LF_{\bbD}^{(\alpha_1, 0), (\beta, 1)}$ and weight its law by $\eps^{\frac12(\alpha_2^2 - \alpha_1^2)} e^{(\alpha_2 - \alpha_1)(\phi, \theta_\eps)}$. Then the law of $\phi|_{\bbD \backslash B_\eps(0)}$ is $\widetilde\LF_{\bbD, \eps}^{(\alpha_2, 0), (\beta, 1)}$.
\end{lemma}

\begin{proof}
Recall $P_\bbD$ is the law of the free boundary GFF on $\bbD$ normalized to have average 0 on $\partial \bbD$. 
    By Girsanov's theorem, for $h$ sampled from $P_\bbD$ weighted by $\eps^{\frac12\alpha^2}e^{\alpha (h, \theta_\eps)}$, we have $h|_{\bbD \backslash B_\eps(0)} \stackrel d= (h' - \alpha \log |\cdot|)|_{\bbD \backslash B_\eps(0)}$ where $h'\sim P_\bbD$.  In other words, this weighting introduces an $\alpha$-log singularity at $0$. Using the above and keeping track of the terms that arise in the definition of the Liouville field, the lemma follows from a direct computation. See \cite[Lemma 4.7]{ARS21} for details in the case where $\alpha_1 = \beta = \gamma$; the argument is identical in our setting. 
\end{proof}

\begin{lemma}\label{lem-weight-fancy}
Let $\alpha_1, \alpha_2, \beta \in \bbR$ and $\eps \in (0,1)$. 
    Let $z \in  \bbD \backslash \{0\}$ and let $K \subset \ol \bbD$ be a compact set such that $\bbD \backslash K$ is simply connected, contains $0$, and has $z$ on its boundary. Let $f: \bbD \backslash K \to \bbD$ be the conformal map such that $f(0) = 0$ and $f(z) = 1$. Let $\widetilde\LF_{\bbD, K, \eps}^{(\alpha_2, 0), (\gamma, z), (\beta, 1)}$ be the law of $\psi|_{\bbD \backslash f^{-1}(B_\eps(0))}$ where $\psi \sim \LF_{\bbD}^{(\alpha_2, 0), (\gamma, z), (\beta, 1)}$.
    \begin{itemize}
        \item Define the {pushforward measure $\hat \theta_\eps = f_*^{-1}\theta_\eps$}. For  $\phi \sim \LF_{\bbD}^{(\alpha_1, 0), (\gamma, z), (\beta, 1)}$ with its law weighted by $({|(f^{-1})'(0)|}\eps)^{\frac12(\alpha_2^2 - \alpha_1^2)} e^{(\alpha_2 - \alpha_1)(\phi, \hat\theta_\eps)}$, the  law of $\phi|_{\bbD \backslash f^{-1}(B_\eps(0))}$ is $\widetilde\LF_{\bbD, K, \eps}^{(\alpha_2, 0), (\gamma, z), (\beta, 1)}$.
    \item Suppose $\alpha_1 + \alpha_2 = 2Q$. For $\phi \sim \LF_\bbD^{(\alpha_1, 0), (\gamma, z), (\beta, 1)}$ with its law weighted by $\eps^{\frac12(\alpha_2^2 - \alpha_1^2)} e^{(\alpha_2 - \alpha_1)(f \bullet_\gamma \phi, \theta_\eps)}$, the law of $\phi|_{\bbD \backslash f^{-1}(B_\eps(0))}$ is $\widetilde\LF_{\bbD, K, \eps}^{(\alpha_2, 0), (\gamma, z), (\beta, 1)}$.
    \end{itemize}
\end{lemma}
\begin{proof}
    The first claim follows from the same argument as that of Lemma~\ref{lem-weight}. Indeed, ${|(f^{-1})'(0)|}\eps$ is the conformal radius of {$f^{-1}(\partial B_\eps(0))$} viewed from $0$ and $\hat \theta_\eps$ is a probability measure on $f(\partial B_\eps(0))$, and  these play the role of $\eps$ and $\theta_\eps$ in Lemma~\ref{lem-weight}. See \cite[Lemma 4.8]{ARS21} for details. For the second claim, note that 
   {\[(f \bullet_\gamma \phi, \theta_\eps) = (\phi \circ f^{-1} + Q \log |(f^{-1})'|, \theta_\eps) = (\phi, f^{-1}_*\theta_\eps) + Q(\log |(f^{-1})'|, \theta_\eps) = (\phi, \hat \theta_\eps) + Q\log |(f^{-1})'(0)|.\]}
    Since $\alpha_1+\alpha_2 = 2Q$ implies $(\alpha_2 - \alpha_1)Q = \frac12(\alpha_2^2 - \alpha_1^2)$, we conclude $({|(f^{-1})'(0)|}\eps)^{\frac12(\alpha_2^2 - \alpha_1^2)} e^{(\alpha_2 - \alpha_1)(\phi, \hat\theta_\eps)} = \eps^{\frac12(\alpha_2^2 - \alpha_1^2)} e^{(\alpha_2 - \alpha_1)(f \bullet_\gamma \phi, \theta_\eps)}$. This with the first claim gives the second claim.
\end{proof}

\begin{figure}
    \centering
    \includegraphics[scale=0.7]{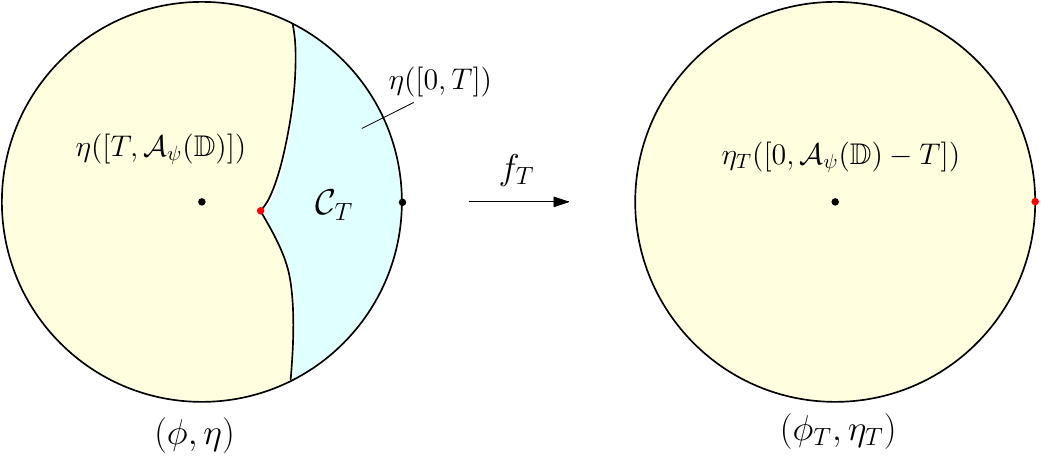}
    \caption{ {An illustration for  Proposition~\ref{prop:radial-finite}. We prove that by cutting the quantum disk $(\bbD,\phi,0,1)$ with the radial $\SLE_\kappa$ curve $\eta$ up to quantum time $T$ restricted to the event $\eta([0,T])$ is simply connected, one gets an independent pair of a quantum cell $\cC_T$ and a quantum disk $(\bbD,\phi_T,0,1)$ after restricting to the event $\{X_T^+(\cC)+Y_T^+(\cC)<\mathcal{L}_{\phi_T}(\partial\bbD)\}$.}}
    \label{fig:radial-finite}
\end{figure}

\begin{proposition}\label{prop:radial-finite}
    Let {$(\phi,\eta, T)$} be a sample from $\mathds{1}_{0<t<\mathcal{A}_{\phi}(\bbD)}\LF_\bbD^{(Q-\frac{\gamma}{4},0),(\frac{3\gamma}{2},1)}(d\phi)\times \mathrm{raSLE}_\kappa(d\eta)\times dt$ and parametrize $\eta$ by its $\mathcal{A}_{\phi}$ quantum area. For $t>0$,  let $f_t:\bbD \backslash \eta([0,t])\to\bbD$ be the conformal map fixing 0 such that $f_t(\eta(t))=1$. Let $\phi_t = f_t\bullet_\gamma\phi$, $\eta_t(s)= f_t(\eta(s + t))$ for $0\le s \le \mathcal{A}_{\phi}(\bbD)-t$, and $\mathcal{C}_t =  (\eta([0,t]), \phi, {\eta}|_{[0,t]})/{\sim_\gamma}$.  Restricted to the event that $\eta([0,T])$ is simply connected, the law of $(\phi_T,\eta_T, \mathcal{C}_T, T)$ is
    \begin{equation}\label{eq:radial-cell-finite}
        {\mathds{1}_{X_t^+(\cC)+Y_t^+(\cC)<\mathcal{L}_{\phi_t}(\partial\bbD)}}\LF_\bbD^{(Q-\frac{\gamma}{4},0),(\frac{3\gamma}{2},1)}(d\phi_t)\times \mathrm{raSLE}_\kappa(d\eta)\times P_t(d\mathcal{C})\,\mathds{1}_{t>0}dt,
    \end{equation}
    where $X_t^+(\cC), Y_t^+(\cC)$ are as in Definition~\ref{def-cell}. 
\end{proposition}

\begin{proof}
 {See Figure~\ref{fig:radial-finite}.}
    Sample $(\widetilde\phi,z,\widetilde\eta)$ from $$\LF_\bbD^{(Q+\frac{\gamma}{4},0),(\gamma,z),(\frac{3\gamma}{2},1)}(d \widetilde \phi)\,\mathrm{raSLE}_\kappa(d\widetilde\eta)\,\mathds{1}_{z\in\bbD}dz$$ 
    and parametrize $\widetilde \eta$ by $\cA_\phi$-quantum area. Let $A$ be the time such that $\widetilde \eta(A) = z$, let  $\widetilde \eta^z = \widetilde \eta|_{[0,A]}$, let $f: \bbD \backslash \widetilde \eta^z \to \bbD$ be the conformal map such that $f(0) = 0$ and $f(z) = 1$, let $\widetilde \phi_A = f \bullet_\gamma \widetilde \phi$, let $\widetilde \eta_A = f \circ \widetilde\eta(\cdot+A)$ and let $\widetilde \cC_A = (\widetilde \eta^z([0,A]), \widetilde \phi, \widetilde \eta^z)/{\sim_\gamma}$. By Lemma~\ref{lem:qtypical} the law of $(\widetilde \phi , \widetilde \eta, A)$ is $\LF_\bbD^{(Q+\frac{\gamma}{4},0),(\frac{3\gamma}{2},1)}(d\widetilde \phi)\times \mathrm{raSLE}_\kappa(d\widetilde\eta)\times \mathds{1}_{a>0} da$ so Proposition~\ref{prop:radial-cell} implies the law of {$(\widetilde\phi_A, \widetilde \eta_A,\widetilde \cC_A, A)$} restricted to the event $\{\widetilde \eta^z \text{ is simply connected}\}$ is 
    \[\mathds{1}_{\widetilde X_a^+(\widetilde \cC) + \widetilde Y_a^+(\widetilde \cC) < \cL_{\widetilde \phi_a}(\partial \bbD)} \LF_{\bbD}^{(Q+\frac{\gamma}{4},0),(\frac{3\gamma}{2},1)}(d\widetilde \phi_a)\times \mathrm{raSLE}_\kappa \times P_a(d\widetilde{\mathcal{C}})\mathds{1}_{a>0}da.\]

       Now for $\e>0$, let $\theta_\e$ be the uniform probability measure on $\partial B_\e(0)$. Let $\alpha_1=Q+\frac{\gamma}{4}$ and $\alpha_2=Q-\frac{\gamma}{4}$. Weight  the law of $(\widetilde \phi,z,\widetilde\eta)$ by $\e^{\frac{\alpha^2_2-\alpha_1^2}{2}}e^{(\alpha_2-\alpha_1)(\widetilde\phi,\theta_\e)}$. 
       By Lemma~\ref{lem-weight-fancy} the law of  $(\widetilde\phi|_{\bbD \backslash \widetilde f_z^{-1}(B_\eps(0))},z,\widetilde\eta)$ under this weighting is 
       \[\widetilde\LF_{\bbD, \widetilde \eta^z, \eps}^{(Q-\frac{\gamma}{4},0),(\gamma,z),(\frac{3\gamma}{2},1)}(d\widetilde \phi)\,\mathrm{raSLE}_\kappa(d\widetilde\eta)\,\mathds{1}_{z\in\bbD}dz.\] 
       On the other hand, by Lemma~\ref{lem-weight}, the weighted law of {$(\widetilde \phi_A|_{\bbD \backslash B_\eps(0)},\widetilde \eta_A, \widetilde \cC_A,  A)$} restricted to the event $\{\widetilde \eta^z \text{ is simply connected}\}$ is
    \eqb\label{eq-cut-some}
    \mathds{1}_{\widetilde X_a^+(\widetilde\cC)+\widetilde Y_a^+(\widetilde\cC)<\mathcal{L}_{\widetilde\phi_a}(\partial\bbD)}\widetilde\LF_{\bbD, \eps}^{(Q-\frac{\gamma}{4},0),(\frac{3\gamma}{2},1)}(d\widetilde\phi_a)\times \mathrm{raSLE}_\kappa \times P_a(d\widetilde{\mathcal{C}})\mathds{1}_{a>0}da.
    \eqe

    To rephrase, if $(\widetilde\phi,z,\widetilde\eta)$ is sampled from
    \begin{equation}\label{eq:law-after-reweight}
        \LF_{\bbD}^{(Q-\frac{\gamma}{4},0),(\gamma,z),(\frac{3\gamma}{2},1)}(d\widetilde\phi)\,\mathrm{raSLE}_\kappa(d\widetilde\eta)\,\mathds{1}_{z\in\bbD}dz
    \end{equation}
    with $\widetilde\eta$ parametrized by quantum area, and $A$ is the time when $\widetilde\eta$ hits $z$, then on the event where $\widetilde \eta^z \text{ is simply connected}$, the law of $(\widetilde\phi_A|_{\bbD \backslash B_\eps(0)},\widetilde \eta_A, \widetilde{\mathcal{C}}_A,A)$ is given by~\eqref{eq-cut-some}. Sending $\eps \to 0$, the same statement holds for $\eps = 0$ when~\eqref{eq-cut-some} is replaced by~\eqref{eq:radial-cell-finite}. 
    On the other hand, by Lemma~\ref{lem:qtypical}, the law of $(\phi,\eta,\eta(T))$ is given by~\eqref{eq:law-after-reweight} (up to renaming of variables). 
     We conclude the proof by observing that the pair $(\phi,\eta,\eta(T))$ and the pair $(\phi,\eta,T)$ uniquely determine each other.      
\end{proof}

Recall the disintegration by quantum boundary length  $\{\LF_{\bbD,\ell}^{(Q-\frac{\gamma}{4},0),(\frac{3\gamma}{2},1)}\}_{\ell>0}$ from Lemma~\ref{lem:disk-boundary-length}.

\begin{corollary}\label{cor:radial-finite}
Fix $t,\ell_0 > 0$. 
    Let $(\phi,\eta)$ be a sample from $\mathds{1}_{\mathcal{A}_{\phi}(\bbD)>t}\LF_{\bbD,\ell_0}^{(Q-\frac{\gamma}{4},0),(\frac{3\gamma}{2},1)}(d\phi)\times \mathrm{raSLE}_\kappa(d\eta)$ and parametrize $\eta$ by its $ {\mathcal{A}_{\phi}}$ quantum area. Let $f_t,\phi_t,\eta_t$ and $\cC_t$ be determined by  $(\phi,\eta)$ in the same way as Proposition~\ref{prop:radial-finite}. 
     Then on the event that $\eta([0,t])$ is simply connected, the law of $(\phi_t,\mathcal{C}_t,\eta_t)$ is
    \begin{equation}
        {\mathds{1}_{X_t^-(\cC)+Y_t^-(\cC)<\ell_0}\LF_{\bbD, \ell_0+X_t(\cC)+Y_t(\cC)}^{(Q-\frac{\gamma}{4},0),(\frac{3\gamma}{2},1)}(d\phi_t)}\,P_t(d\mathcal{C})\times\mathrm{raSLE}_\kappa(d\eta),
    \end{equation}
    where $X_t(\cC) ,Y_t(\cC),X_t^-(\cC), Y_t^-(\cC)$ are as in Definition~\ref{def-cell}. 
\end{corollary}
\begin{proof}
    If we do not fix the boundary length of $\phi$, i.e., we instead assume that $(\phi,\eta)$ is sampled from $\mathds{1}_{\mathcal{A}_{\phi}(\bbD)>t}\LF_{\bbD}^{(Q-\frac{\gamma}{4},0),(\frac{3\gamma}{2},1)}(d\phi)\times \mathrm{raSLE}_\kappa(d\eta)$, then it follows from Proposition~\ref{prop:radial-finite} by disintegrating on the value of $T$ that the law of $(\phi_t,\mathcal{C}_t,\eta_t)$ is
    \begin{equation}
         {\mathds{1}_{X_t^+(\cC)+Y_t^+(\cC)<\mathcal{L}_{\phi_t}(\partial\bbD)}}\LF_{\bbD}^{(Q-\frac{\gamma}{4},0),(\frac{3\gamma}{2},1)}(d\phi_t)\,P_t(d\mathcal{C})\times\mathrm{raSLE}_\kappa(d\eta).
    \end{equation}
    Now we disintegrate over $\mathcal{L}_{\phi}(\partial\bbD)$, and the claim follows from $\mathcal{L}_{\phi_t}(\partial\bbD) = \mathcal{L}_{\phi}(\partial\bbD) + {X_t(\cC) + Y_t(\cC)}$ {and $\{ X_t^+(\cC)+Y_t^+(\cC)<\mathcal{L}_{\phi_t}(\partial\bbD)\}  = \{X_t^- (\cC) + Y_t^- (\cC) < \mathcal L_{\phi}(\partial \bbD) \}$}.
\end{proof}

\begin{proof}[Proof of Theorem~\ref{thm-mot-disk-finite}]
Recall that $\mathrm{BM}_\kappa$ is the law of  correlated two-dimensional Brownian motion $(\widetilde X_t, \widetilde Y_t)_{t \geq 0}$ with $\widetilde X_0 = \widetilde Y_0=0$ and covariance given by~\eqref{eqn-thm:mot-chordal}. Sample {$(\widetilde X_t,\widetilde Y_t)_{t \geq 0}$} from $\mathrm{BM}_\kappa$ and let 
     $\widetilde \tau$ be the first time $t$ that $1 + \widetilde X_t + \widetilde Y_t = 0$. Our first goal is to show that $(X_t, Y_t)_{[0,A]} \stackrel d= (\widetilde X_t, \widetilde Y_t)_{[0, \widetilde \tau]}$. To that end, we will show that $(X_s, Y_s)_{[0, \tau_1]} \stackrel d= (\widetilde X_s, \widetilde Y_s)_{[0,\widetilde \tau_1]}$ for suitable stopping times $\tau_1, \widetilde \tau_1$ corresponding to ``wrapping around'', then iterate to conclude. Afterwards, we establish~\eqref{eq-thm-radial} to complete the proof.

{Recall that $Z := |\LF_{\bbD, \ell}^{(Q - \frac\gamma4, 0), (\frac{3\gamma}2, 1)}|$ does not depend on $\ell$ (Lemma~\ref{lem:disk-boundary-length}).}
    Suppose $\phi$ is a sample from  $Z^{-1}\LF_{\bbD,1}^{(Q-\frac{\gamma}{4},0),(\frac{3\gamma}{2},1)}$, and $\eta$ is an independent radial $\SLE_\kappa$ process from 1 to 0 parametrized by its $\mathcal{A}_\phi$ quantum area. {Fix $t > 0$ and} let $f_t,\phi_t,\eta_t$ and $\cC_t$ be determined by  $(\phi,\eta)$ in the same way as Proposition~\ref{prop:radial-finite}. 
    By Corollary~\ref{cor:radial-finite}, when restricted to the {event $F_t$ that $\cA_{\phi}(\bbD)>t$ and} $\eta([0,t])$ is simply connected, the joint law  of $(\phi_t,\mathcal{C}_t,\eta_t)$ is
    \begin{equation*}
        {\mathds{1}_{X_t^-(\cC)+Y_t^-(\cC)<1}Z^{-1}\LF_{\bbD, 1+X_t(\cC)+Y_t(\cC)}^{(Q-\frac{\gamma}{4},0),(\frac{3\gamma}{2},1)}(d\phi_t)}\,P_t(d\mathcal{C}_t)\times\mathrm{raSLE}_\kappa(d\eta_t), 
    \end{equation*}
     {so the joint law of $(\phi_t,{(X_\cdot,Y_\cdot)|_{[0,t]}},\eta_t)$ is}
     \begin{equation*}
        Z^{-1}\LF_{\bbD, 1+X_t+Y_t}^{(Q-\frac{\gamma}{4},0),(\frac{3\gamma}{2},1)}(d\phi_t)\,\times {\mathds{1}_{\widetilde F_t}\mathrm{BM}_\kappa^t(d(X_\cdot,Y_\cdot))}\times\mathrm{raSLE}_\kappa(d\eta_t)
    \end{equation*}
    {where $\mathrm{BM}_\kappa^t$ is the law of a sample from $\mathrm{BM}_\kappa$ restricted to the time interval $[0,t]$, and $\widetilde F_t= \{ -\inf_{[0,t]} X_\cdot - \inf_{[0,t]} Y_\cdot > -1\}$. Since $|Z^{-1}\LF_{\bbD, 1+X_t+Y_t}^{(Q-\frac{\gamma}{4},0),(\frac{3\gamma}{2},1)}| = |\mathrm{raSLE}_\kappa| = 1$ regardless of the value of $1+X_t + Y_t$, the marginal law of $(X_\cdot, Y_\cdot)|_{[0,t]}$ restricted to $F_t$ is $\mathds{1}_{\widetilde F_t} \mathrm{BM}_\kappa^t$. Since $t$ is arbitrary, we conclude that $(X_\cdot, Y_\cdot)_{[0,\tau_1]} \stackrel d= (\widetilde X_\cdot, \widetilde Y_\cdot)_{[0,\widetilde \tau_1]}$ where $\tau_1 = \inf\{ s : -\inf_{[0,s]} X_\cdot - \inf_{[0,s]}Y_\cdot \leq -1\}$ and $\widetilde \tau_1 = \inf\{ s : -\inf_{[0,s]} \widetilde X_\cdot - \inf_{[0,s]}\widetilde Y_\cdot \leq -1\}$.} 

{Next, let  $\tau_2$ (resp.\ $\widetilde\tau_2$) be the first time $t>\tau_1$ (resp.\ $t>\widetilde\tau_1$) that $\inf_{\tau_1<s<t}X_s+\inf_{\tau_1<s<t}Y_s = -1$ (resp.\ $\inf_{\widetilde\tau_1<s<t}\widetilde X_s+\inf_{\widetilde\tau_1<s<t}\widetilde Y_s = -1$). We will show that $(X_\cdot, Y_\cdot)_{[0,\tau_2]} \stackrel d= (\widetilde X_\cdot, \widetilde Y_\cdot)_{[0,\widetilde\tau_2]}$.}
    Fix $t_1>0$ and condition on $\{t_1<\tau_1\}$. Then the conditional law of $(\phi_{t_1}, \eta_{t_1})$ given $\mathcal{C}_{t_1}$ is ${Z^{-1}}\LF_{\bbD,1+X_{t_1}+Y_{t_1}}^{(Q-\frac{\gamma}{4},0),(\frac{3\gamma}{2},1)}(d\phi)\times \mathrm{raSLE}_\kappa$, and the boundary length process of $(\phi_{t_1}, {\eta_{t_1}})$ is specified by $(X_t-X_{t_1},Y_t-Y_{t_1})_{t_1\le t\le A}$. 
    Therefore following the same reasoning, if we let $\sigma_2$ (resp.\ $\widetilde\sigma_2$) be the first time $t$ such that $\inf_{t_1<s<t}(X_s-X_{t_1})+\inf_{t_1<s<t}(Y_s-Y_{t_1})=-1-X_{t_1}-Y_{t_1}$ (resp.\  $\inf_{t_1<s<t}(\widetilde X_s-\widetilde X_{t_1})+\inf_{t_1<s<t}(\widetilde Y_s-\widetilde Y_{t_1})=-1-\widetilde X_{t_1}-\widetilde Y_{t_1}$), then $(X_t-X_{t_1},Y_t-Y_{t_1})_{t_1\le t\le \sigma_2}$ is independent of $(X_s,Y_s)_{0\le s\le t_1}$ and agrees in law with $(\widetilde X_t-\widetilde X_{t_1},\widetilde Y_t-\widetilde Y_{t_1})_{t_1\le t\le \widetilde\sigma_2}$ conditioned on $\{t_1<\widetilde\tau_1\}$. This implies that conditioned on $\{t_1<\tau_1\}$, the law of $(X_s,Y_s)_{0\le s\le \sigma_2}$ agrees with that of $(\widetilde X_s,\widetilde Y_s)_{0\le s\le \widetilde\sigma_2}$  conditioned on $\{t_1<\widetilde\tau_1\}$. Since $t_1$ is  arbitrary, we conclude  $(X_\cdot,Y_\cdot)_{[0,\tau_2]} \stackrel d= (\widetilde X_\cdot,\widetilde Y_\cdot)_{[0,\widetilde\tau_2]}$. 

   {Arguing similarly, if we iteratively define $\tau_n$ (resp.\ $\widetilde \tau_n$) to be the first time $t > \tau_{n-1}$ (resp.\ $t > \widetilde \tau_{n-1}$) such that $\inf_{\tau_{n-1} < s  <t} X_s + \inf_{\tau_{n-1} < s < t} Y_s = -1$ (resp.\ $\inf_{\widetilde\tau_{n-1} < s  <t} \widetilde X_s + \inf_{\widetilde\tau_{n-1} < s < t} \widetilde Y_s = -1$), then $(X_\cdot, Y_\cdot)_{[0,\tau_n]} \stackrel d= (\widetilde X_\cdot, \widetilde Y_\cdot)_{[0,\widetilde \tau_n]}$ for all $n$. Since $\lim_{n \to \infty} \tau_n = A$ and $\lim_{n \to \infty} \widetilde \tau_n = \widetilde \tau$} where $\widetilde\tau=\inf\{t>0:1+\widetilde X_t+\widetilde Y_t=0\}$, it follows that $(X_t,Y_t)_{0\le t\le A}\stackrel d= (\widetilde X_t,\widetilde Y_t)_{0\le t\le \widetilde\tau}$. This proves the first claim.

Finally, we prove~\eqref{eq-thm-radial}, 
{which is immediate from Corollary~\ref{cor:radial-finite} when $s=0$. We first claim that for each fixed $s>0$, conditioned on the event $s<\cA_\phi(\bbD)$ and $(X_\cdot,Y_\cdot)|_{[0,s]}$, the law of 
$(\phi_s, \eta_s)$ is
$Z^{-1}\LF_{\bbD,1+X_s+Y_s}^{(Q - \frac\gamma4, 0), (\frac{3\gamma}2, 1)} \times \mathrm{raSLE}_\kappa$. To see this, fix $n>0$. For $1\le k\le 2^n$, let $E_{n,k,s}$ be the event where $\frac{ks}{2^n}<\cA_\phi(\bbD)$ and for each $1\le j\le k$, $\eta([\frac{(j-1)s}{2^n},\frac{js}{2^n}])$ is simply connected. Then conditioned on $E_{n,1,s}$ and $(X_\cdot,Y_\cdot)|_{[0,\frac{s}{2^n}]}$, by Corollary~\ref{cor:radial-finite} the law of $(\phi_{\frac{s}{2^n} }, \eta_{\frac{s}{2^n}})$ is $Z^{-1}\LF_{\bbD,1+X_{\frac{s}{2^n}}+Y_{\frac{s}{2^n}}}^{(Q - \frac\gamma4, 0), (\frac{3\gamma}2, 1)} \times \mathrm{raSLE}_\kappa$. Applying Corollary~\ref{cor:radial-finite} once more to $(\phi_{\frac{s}{2^n} }, \eta_{\frac{s}{2^n}})$, we see that conditioned on $E_{n,2,s}$ and $(X_\cdot,Y_\cdot)|_{[0,\frac{2s}{2^n}]}$,  the law of $(\phi_{\frac{2s}{2^n} }, \eta_{\frac{2s}{2^n}})$ is $Z^{-1}\LF_{\bbD,1+X_{\frac{2s}{2^n}}+Y_{\frac{2s}{2^n}}}^{(Q - \frac\gamma4, 0), (\frac{3\gamma}2, 1)} \times \mathrm{raSLE}_\kappa$. By  iterating this argument $2^n$ times,  conditioned on $E_{n,2^n,s}$, the law of 
$(\phi_s, \eta_s)$ is
$Z^{-1}\LF_{\bbD,1+X_s+Y_s}^{(Q - \frac\gamma4, 0), (\frac{3\gamma}2, 1)} \times \mathrm{raSLE}_\kappa$. On the other hand, using the continuity of the curve $\eta$,  conditioned on $\cA_{\phi}(\bbD)>s$ the event $E_{n,2^n,s}$ holds with probability $1-o_1(n)$ as $n\to\infty$. Now we can apply  Corollary~\ref{cor:radial-finite} to $(\phi_s,\eta_s)$ and }conclude that conditioned on the event that $t<\cA_{\phi}(\bbD)$ and $\eta_s([0,t-s])$ is simply connected, the law of $(\eta_s([0,t-s]), \phi_s, \eta_s|_{[0,t-s]})/{\sim_\gamma}$ is absolutely continuous with respect to $P_{t-s}$. Therefore  $F((X_{\cdot + s}-X_{s}, Y_{\cdot + s}-Y_{s})|_{[0, t-s]}) = (\eta_s([0,t-s]), \phi_s, \eta_s|_{[0,t-s]})/{\sim_\gamma}$ a.s.. By definition $(\eta([s,t]), \phi, \eta(\cdot +s)|_{[0,t-s]})/{\sim_\gamma} = (\eta_s([0,t-s]), \phi_s, \eta_s|_{[0,t-s]})/{\sim_\gamma}$, so~\eqref{eq-thm-radial} holds.

\end{proof}

\section{A spherical mating-of-trees}\label{sec-main-proof}
In this section we prove Theorem~\ref{thm-main}. The key ingredient is the following spherical mating-of-trees result of independent interest. 
Recall $\mathcal M_2^\mathrm{sph}(\alpha)$ is the law of the quantum sphere from Definition~\ref{def-sphere}.

Suppose $\kappa \geq 8$ and $\gamma = \frac4{\sqrt\kappa}$. 
Let $(\bbC, \phi, \infty,0)$ be an embedding of a sample from $\mathcal M_2^\mathrm{sph}(Q - \frac\gamma4)$ conditioned to have quantum area at least 1; write $A = \cA_\phi(\bbC)$. Let $\eta:[0, A]\to \hat\bbC$ be an independent whole-plane $\SLE_{\kappa}$ from $0$ to $\infty$ parametrized by quantum area. There is a unique continuous process $(X_t, Y_t)_{[0,A]}$ starting at $(X_0, Y_0) = (0,0)$ which keeps track of the changes in the left and right boundary lengths of $\eta([0,t])$, in the following sense. For any $s \in (0,A)$ and point $p \in \partial (\eta([0,s]))$ different from $\eta(s)$, let $\sigma>s$ be the next time $\eta$ hits $p$. For each $t \in (s, \sigma)$ let $X_t^s$ (resp.\ $Y_t^s$) be the quantum length of the counterclockwise (resp.\ clockwise) boundary arc of $\eta([0,s])$ from $\eta(s)$ to $p$. Then $(X_t^s - X_s^s, Y_t^s - Y_s^s)_{(s,\sigma)} = (X_t - X_s, Y_t - Y_s)_{(s, \sigma)}$.  See Figure~\ref{fig:bdrylength-sph} for an illustration. To justify the existence and uniqueness of $(X_t, Y_t)_{[0,A]}$, similarly to the radial case in Theorem~\ref{thm-mot-disk-finite} we can define a process $(\widetilde X_t, \widetilde Y_t)_{(0,A)}$ with this property and which is unique up to additive constant. We extend it to $(\widetilde X_t, \widetilde Y_t)_{[0,A]}$ by continuity, and thus  uniquely fix $(X_t, Y_t)_{[0,A]} = (\widetilde X_t - \widetilde X_0, \widetilde Y_t - \widetilde Y_0)_{[0,A]}$.
\begin{theorem}[Spherical mating-of-trees]\label{thm-sphere-mot}
 Let $(L_t, Z_t) = (X_t + Y_t, X_t - Y_t)$. Then $L_t$ has the law of a Brownian excursion with quadratic variation $(2\mathbbm a \sin(\frac{\pi\gamma^2}{8}))^2\, dt$ conditioned to have duration at least 1, and given the process $(L_t)$ with random duration $\tau$, the process $(Z_t)_{[0,\tau]}$ is conditionally independent Brownian motion with quadratic variation $(2\mathbbm a \cos(\frac{\pi\gamma^2}8))^2\,dt$ run for time $\tau$. Here $\mathbbm a$ is as in~\eqref{eqn-thm:mot-chordal}. 
Moreover, for any $0 < s<t$, on the event that $t <\tau$ and $\eta([s,t])$ is simply connected, we have 
\eqb\label{eq-thm-sph}
F((X_{\cdot + s} - X_s, Y_{\cdot + s} - Y_s)_{[0,t-s]}) = (\eta([s,t]), \phi, \eta(\cdot + s)|_{[0,t-s]})/{\sim_\gamma}
\eqe
where $F$ is the map from Lemma~\ref{lem-F-reversible}.
\end{theorem}

\begin{figure}
    \centering
    \includegraphics[scale=0.4]{figures/lengths-sphere.pdf}
    \caption{The boundary length process $(X_t, Y_t)_{[0,A]}$ of Theorem~\ref{thm-sphere-mot} is characterized by $X_0 = Y_0 = 0$ and the property that for each time $s$ and choice of boundary point $p \in \partial \eta([0,s])$ not equal to $\eta(s)$, for any time $t > s$ before  the time $\eta$ next hits $p$, we have $(X_t^s - X_s^s, Y_t^s - Y_s^s) = (X_t - X_s, Y_t - Y_s)$, where $(X_\cdot^s, Y_\cdot^s)$ is shown in red and blue. Here $\eta([0,s])$ is shown in dark gray, and $\eta([s,t])$ is colored light grey.}
    \label{fig:bdrylength-sph}
\end{figure}
We note that $L_t$ is the quantum length of $\partial (\eta([0,t]))$ for all $t$. 

To prove Theorem~\ref{thm-sphere-mot}, we start with the radial mating of trees Theorem~\ref{thm-mot-disk-finite}, and condition on having quantum area at least 1 but having small boundary length $\ell\ll 1$ (event $F_\ell$ from~\eqref{eq-F}). 
Lemma~\ref{lem-limit-BM} below implies that when $\ell \to 0$ the  limiting boundary length process is that of Theorem~\ref{thm-sphere-mot}. On the other hand, when $\ell \to 0$ the curve-decorated quantum surface converges to the conditioned quantum sphere decorated by independent whole-plane SLE (Proposition~\ref{prop-sphere-given-F}). Combining these two facts gives Theorem~\ref{thm-sphere-mot}.

\begin{lemma}\label{lem-limit-BM}
Let $\ell > 0$ and let $L_t^\ell$ be Brownian motion starting at $\ell$ and having quadratic variation $(2\mathbbm a \sin(\frac{\pi\gamma^2}{8}))^2\, dt$, run until the time $\tau$ that it first hits 0. Given $(L_t^\ell)_{[0,\tau]}$ let $Z_t^\ell$ be an independent Brownian motion with quadratic variation $(2\mathbbm a \cos(\frac{\pi\gamma^2}{8}))^2\, dt$ run for time $\tau$. 
As $\ell \to 0$, the process $(L_t^\ell, Z_t^\ell)$ conditioned on $\tau\ge1$ converges in distribution to the Brownian process described in Theorem~\ref{thm-sphere-mot}.
\end{lemma}
\begin{proof}
This is immediate from the limiting construction of the Brownian excursion. 
\end{proof}

Given Theorem~\ref{thm-sphere-mot}, the proof of Theorem~\ref{thm-main} goes as follows. Let $\mathfrak S$ be the SLE-decorated quantum surface in Theorem~\ref{thm-sphere-mot}. Since its boundary length process agrees in law with its time-reversal, we have $\mathfrak S \stackrel d= \widetilde {\mathfrak S}$ where $\widetilde{\mathfrak S}$ is obtained from $\mathfrak S$ by switching its two points and reversing its curve. This implies that the law of the curve is reversible, as desired. 

In Section~\ref{subsec-pinch} we show that a certain quantum sphere can be obtained from a disk by taking a limit (Proposition~\ref{prop-sphere-given-F}). In Section~\ref{subsec-proofs-main} we use this to obtain Theorem~\ref{thm-sphere-mot} and then Theorem~\ref{thm-main}.

\subsection{Pinching an LCFT disk to get an LCFT sphere}\label{subsec-pinch}

The goal of this section is to prove Proposition~\ref{prop-sphere-given-F} which states that a Liouville field on the disk conditioned to have area at least 1 and boundary length $\ell$ converges as $\ell \to 0$ to {a sample from $\mathcal M_2^\mathrm{sph}(Q -\frac\gamma4)$} conditioned to have  area at least 1. Although the statement of Proposition~\ref{prop-sphere-given-F} does not involve SLE or mating-of-trees, our arguments will use these to establish that the field remains ``well behaved'' near the boundary despite the conditioning on low probability events. 

Instead of working in the domains $\bbC$ and $\bbD$, we will parametrize by the horizontal cylinder $\cC := (\bbR \times [0,2\pi])/{\sim}$ and half-cylinder $\cC_+ := ([0, \infty) \times [0,2\pi])/{\sim}$ where the upper and lower boundaries are identified by $x\sim x + 2\pi i$. This simplifies our exposition later. 

We first define the Liouville field on $\cC_+$.  Let $f:\cC_+ \to \bbD$ be the map such that $f(z) = e^{-z}$.
\begin{definition}\label{def-LF-C}
For $\alpha, \beta \in \bbR$ and $\ell >0$, define $\LF_{\cC_+, \ell}^{(\alpha, +\infty), (\beta, 0)} := f^{-1} \bullet_\gamma \LF_{\bbD, \ell}^{(\alpha, 0), (\beta, 1)}$. 
\end{definition}

 $\LF_{\cC_+, \ell}^{(\alpha, +\infty), (\beta, 0)}$ inherits the following Markov property from $\LF_{\bbD, \ell}^{(\alpha, 0), (\beta, 1)}$. For $\phi \sim \LF_{\cC_+, \ell}^{(\alpha, +\infty), (\beta, 0)}$, conditioned on $\phi|_{\partial \cC_+}$ we have
\eqb\label{eq-markov-cC}
\phi \stackrel d= \mathfrak h + h_0 -(Q-\alpha) \Re \cdot,
\eqe
where $\mathfrak h$ is the harmonic function on $\cC_+$ with boundary conditions $\phi|_{\cC_+}$, and $h_0$ is a Dirichlet GFF on $\cC_+$.

We define a probability measure $\cL$ on fields on $\mathcal C$ as follows. 
Consider $(\hat h, \mathbf c)$ sampled as in Definition~\ref{def-sphere} with $\alpha = Q - \frac\gamma4$ and conditioned on the event that $\cA_{\hat h + \mathbf c} (\cC) > 1$, let $\sigma \in \bbR$ satisfy $\cA_{\hat h + \mathbf c} ([\sigma, +\infty) \times [0,2\pi]) = \frac12$,  let $\phi' = \hat h (\cdot +\sigma) + \mathbf c$, and let $\cL$ be the law of $\phi'$. 
Thus,  $\phi'\sim \cL$ corresponds to a sample from $\mathcal M_2^\mathrm{sph}(\alpha)$ conditioned to have quantum area greater than 1, embedded such that $\cA_\phi(\mathcal C_+) = \frac12$. 

The main result of this section is that for small $\ell$,  a field sampled from $\LF_{\mathcal C_+, \ell}^{(\alpha, +\infty), (\beta, 0)}$ conditioned on $F_\ell$ resembles a quantum sphere conditioned to have quantum area at least 1.
\begin{proposition}
\label{prop-sphere-given-F}
Let $(\alpha, \beta) = (Q - \frac\gamma4, \frac{3\gamma}2)$ and $\ell > 0$. Sample $\phi$ from $\LF_{\cC_+, \ell}^{(\alpha, +\infty), (\beta, 0)}$ conditioned on 
    \eqb
\label{eq-F}
	F_{\ell} := \{ \cA_\phi(\cC_+) > 1 \}.
    \eqe
    Let $\sigma>0$ satisfy $\cA_{\phi}(\cC_+ + \sigma) = \frac12$ and let $\widetilde \phi = \phi(\cdot + \sigma)$.  For any $U \subset \mathcal C$ bounded away from $-\infty$, as $\ell \to 0$ the field $\widetilde \phi|_U$ converges in distribution to  $\phi'|_U$ where $\phi' \sim \cL$.
\end{proposition}

{
We first state a version of Proposition~\ref{prop-sphere-given-F} where we additionally condition on the field near $\partial \cC_+$ not behaving too wildly, in the sense that it has ``scale $\ell$'' observables near $\partial \cC_+$. 
}


\begin{lemma}\label{lem-sphere-given-E}
Let $(\alpha, \beta) = (Q - \frac\gamma4, \frac{3\gamma}2)$. Fix a nonnegative smooth function $\rho$ in $\mathcal C$ supported on $[1,2] \times [0,\pi]$, such that $\rho$ is constant on each vertical segment\footnote{This is convenient for the proof of Lemma~\ref{lem-sphere-given-E} since $(\phi, \rho)$ only depends on the projection of $\phi$ to $H_\mathrm{av}(\mathcal C)$.} $\{t\} \times [0,\pi]$ and $\int \rho = 1$.
Let $K, \ell > 0$.
Sample a field $\phi$ from $\LF_{\cC_+, \ell}^{(\alpha, +\infty), (\beta, 0)}$ conditioned on 
\eqb
	\label{eq-E}
	E_{\ell,K} := F_\ell \cap \{ \cA_{\phi - \frac2\gamma \log \ell}([0,1] \times [0,2\pi]) < K \text{ and } |(\phi, \rho)-\frac2\gamma \log \ell| < K\}.
\eqe
Let $\sigma>0$ satisfy $\cA_{\phi}(\cC_+ + \sigma) = \frac12$ and let $\widetilde \phi = \phi(\cdot + \sigma)$.  For any $U \subset \mathcal C$ bounded away from $-\infty$, as $\ell \to 0$ the field $\widetilde \phi|_U$ converges in distribution to  $\phi'|_U$ where $\phi' \sim \cL$.
\end{lemma}


The statement of Lemma~\ref{lem-sphere-given-E} is parallel to that of \cite[Proposition 4.1]{MS19}, except that we condition on an event measurable with respect to $\phi|_{[0,2]\times[0,2\pi]}$ (the second set in RHS of~\eqref{eq-E}), while they more strongly assert the asymptotic independence of $\phi|_{[0,2]\times[0,2\pi]}$ and $\widetilde \phi|_U$ (or rather, the corresponding fields in their setting). 
Using the Markov property~\eqref{eq-markov-cC} of $\LF_{\cC_+, \ell}^{(\alpha, +\infty), (\beta, 0)}$, the proof of Lemma~\ref{lem-sphere-given-E} is identical to the proof of \cite[Proposition 4.1]{MS19}, so we omit it.

Next, we show that conditioned on $F_\ell$, with high probability $E_{\ell, K}$ occurs. To that end, we will control the field near $\partial \cC_+$ when we condition on $F_\ell$ by using the following lemma. Any planar domain $A$ with the annulus topology is conformally equivalent to $\{z: 1 < |z| < e^{2\pi M}\}$ for some unique $M>0$; this $M$  is called the \emph{modulus} of $A$, and we denote it by $\mathrm{Mod}(A)$. 

\begin{lemma}\label{lem-thick-ann}
Let $(\alpha, \beta) = (Q - \frac\gamma4, \frac{3\gamma}2)$ and $n\geq1$. Consider the setting of Theorem~\ref{thm-mot-disk-finite}, except we embed in $(\cC_+, +\infty, 0)$ rather than $(\bbD, 0, 1)$, so $\phi$ is sampled from $\LF_{\cC_+, 1}^{(\alpha, +\infty), (\beta, 0)}$ and $\eta$ is an independent radial $\SLE$ in $(\cC_+, +\infty, 0)$. Let $(L_t, Z_t) = (1+ X_t + Y_t, X_t - Y_t)$ and let $\tau_x = \inf\{ t: L_t = x\}$. Conditioned on $\{ \tau_{2^n} < \tau_{0}\}$, the explored region $A = \eta([0,\tau_{2^n}])$ is annular with probability $1-o_n(1)$, and its modulus tends to $\infty$ in probability as $n \to \infty$.
\end{lemma}
\begin{proof}
First consider $n=1$.
Condition on $\tau_2 < \tau_{0}$ and let $A_1 = \eta([0,\tau_2])$. 
Since Brownian motion stays arbitrarily close to any deterministic path with positive probability, $A_1$ is annular with positive probability. Thus  there exists $m_0>0$ such that the event $E_1 = \{A_1 \text{ annular and } \mathrm{Mod}(A_1) > m_0\}$ has conditional probability $p>0$ {given $\tau_2 < \tau_0$}. 

Now consider general $n \geq 1$. Condition on $\tau_{2^n} < \tau_{0}$, and for $1\leq i \leq n$ define $A_i = \eta([\tau_{2^{i-1}},\tau_{2^i}])$ and $E_i = \{A_i \text{ annular and } \mathrm{Mod}(A_i) > m_0\}$. By the scale invariance and strong Markov property of Brownian motion, the events $E_1, \dots, E_n$ are {conditionally} independent and each occur with probability $p$. Let $I$ be the random set of $i$ such that $E_i$ holds, then $|I| \to \infty$ in probability as $n \to \infty$, so in particular $\bbP[A \text{ annular}] \geq \bbP[|I| \geq 1] \to 1$ as $n \to \infty$. Finally, by the subadditivity of moduli, on the event $\{A \text{ annular}\}$ we have 
\[\mathrm{Mod}(A) \geq \sum_{i\in I} \mathrm{Mod}(A_i) \geq m_0|I|.\]
Since $|I| \to \infty$ in probability, $\mathrm{Mod}(A) \to \infty$ in probability as desired. 
\end{proof}

Lemma~\ref{lem-thick-ann}
 states that on the rare event that the boundary length hits $2^n$, with high probability the explored region $A$ at this hitting time is an annulus with large modulus. Next, we give a uniform bound on the field for all embeddings of $(A, \phi,0)/{\sim}$ in $\cC_+$ having $\partial \cC_+$ as a boundary component.

\begin{lemma}\label{lem-distortion}
There is an absolute constant $m>0$ such that the following holds. 
Fix $n \geq 1$ and let $\rho$ be as defined {as in Lemma~\ref{lem-sphere-given-E}}. In the setting of Lemma~\ref{lem-thick-ann}, condition on $\{ \tau_{2^n } < \tau_{0}\}$ and on $\mathrm{Mod}(A) > m$. Let $\widetilde A \subset \cC_+$ be any bounded annulus having $\partial \cC_+$ as a boundary component such that $\mathrm{Mod}(\widetilde A) = \mathrm{Mod}(A)$. Let $\widetilde \phi$ be the field on $\widetilde A$ such that $(\widetilde A, \widetilde \phi, 0)/{\sim_\gamma}=(A, \phi, 0)/{\sim_\gamma}$, then 
\eqb\label{eq-sup-field}
\sup |(\widetilde \phi, \rho)| < \infty \quad \text{ almost surely}, 
\eqe
where the supremum is taken over all choices of $\widetilde A$. 
\end{lemma}
\begin{proof}
We first fix a canonical embedding $(\widetilde A_0, \widetilde \phi_0, 0)$ by specifying that $A_0$ is concentric, i.e., $A_0 = [0,t] \times [0,2\pi]$  where $t = 2\pi \mathrm{Mod}(A)  > 2 \pi m$. For $b>0$, let $H_b$ be the set of nonnegative smooth functions $f$ in $\mathcal C$ supported in $[\frac12, 3] \times [0,2\pi]$ with $f \geq 0$, $\int f(x) \,dx = 1$ and $\| f'\|_\infty \leq b$. 
Since $\widetilde \phi_0$ is locally absolutely continuous with respect to a GFF, as explained in the paragraph just after \cite[Proposition 9.19]{DMS14} we have  $\sup_{f \in H_b} |(\widetilde \phi_0, f)|<\infty$ almost surely. 

For any other embedding $(\widetilde A, \widetilde \phi, 0)$ let $g: \widetilde A_0 \to \widetilde A$ be the conformal map such that $\widetilde \phi = g \bullet_\gamma \widetilde \phi_0$, then 
\eqb\nonumber
(\widetilde \phi, \rho) = (\widetilde \phi_0 \circ g^{-1} + Q \log |(g^{-1})'|,\rho) = (\widetilde \phi_0, |g'|^2 \rho \circ g) + Q (\log | (g^{-1})'|, \rho). \eqe
Assuming the absolute constant $m$ is chosen sufficiently large, conformal distortion estimates (e.g.\ 
\cite[Theorem 5]{annular-distortion}) give $\sup_{[\frac12,3]\times [0,2\pi]} |g'-1 | < \frac1{10}$. Thus, $|(\widetilde\phi_0, |g'|^2 \rho \circ g)| \leq \sup_{f \in F_{b}}|(\widetilde \phi_0, f)|$ for some $b$ depending only on $\rho$ and $| Q (\log |(g^{-1})'|, \rho)| \leq 10Q$, giving the desired uniform bound for $|(\widetilde \phi, \rho)|$.
\end{proof}

Now, we will prove that conditioned on $F_\ell$, the event $E_{\ell,K}$ is likely. Briefly, conditioning on $F_\ell$, Theorem~\ref{thm-mot-disk-finite} gives a description of the quantum surface near $\partial \cC_+$ which we use to  bound the field average near $\partial \cC_+$ via Lemma~\ref{lem-distortion}.
\begin{proposition}\label{prop-F|E}
Let $(\alpha, \beta) = (Q - \frac\gamma4, \frac{3\gamma}2)$. For each $\delta > 0$ there exists $K_0 > 0$ such that for all $K > K_0$ 
\[\liminf_{\ell \to 0} \LF_{\cC_+, \ell}^{(\alpha, +\infty), (\beta, 0)}[E_{\ell,K} \mid F_{\ell} ] > 1-\delta.\]
\end{proposition}
\begin{proof}

Fix $n = n(\delta) \geq 1$ sufficiently large  such that in the setting of Lemma~\ref{lem-thick-ann} we have $\bbP[\mathrm{Mod}(A) \geq m] \geq 1-\frac{\delta}{4}$,  where $m$ is the absolute constant in Lemma~\ref{lem-distortion}.

Sample $\phi \sim \LF_{\cC_+, \ell}^{(\alpha, +\infty), (\beta, 0)}$, and independently sample a radial $\SLE_\kappa$ curve $\eta$ in $(\cC_+, 0)$ targeting $+\infty$ and parametrized by $\cA_\phi$. The law of  $\phi^0 := \phi - \frac2\gamma \log \ell$ is $\LF_{\cC_+,1}^{(\alpha, +\infty), (\beta, 0)}$. Let  $(X_t, Y_t)$ be the boundary length process for $(\phi^0, \eta)$ as in Lemma~\ref{lem-thick-ann}, let $(L_t, Z_t)  = (1 + X_t + Y_t, X_t - Y_t)$, and let $\tau_x$ be the time $L_t$ first hits $x$. By Lemma~\ref{lem-limit-BM} we have $\LF_{\cC_+, \ell}^{(\alpha, +\infty), (\beta, 0)}[ \tau_{2^n} < \tau_0 \mid F_\ell] = 1-o_\ell(1)$, and furthermore conditioning on $\{ \tau_{2^n} < \tau_0 \} \cap F_\ell$ the conditional  law of $(L_t, Z_t)_{[0,\tau_{2^n}]}$ is within $1-o_\ell(1)$ in total variation distance of the corresponding process of Lemma~\ref{lem-thick-ann}. We conclude that conditioned on $F_\ell$, the conditional law of the quantum surface $\cA := (\eta([0,\ell^{2} \tau_{2^n}]), \phi - \frac2\gamma \log \ell, 0)$ is within $1-o_\ell(1)$ in total variation distance of the quantum surface of Lemma~\ref{lem-thick-ann}, and hence within $1 - \frac\delta2 - o_\ell(1)$ in total variation distance of the quantum surface of Lemma~\ref{lem-distortion}. 
Choose $K_0$ sufficiently large that 
in Lemma~\ref{lem-distortion} the finite constant in~\eqref{eq-sup-field} is bounded by $K_0 - \log 2$ with probability at least $1-\frac\delta4$, and the quantum area of the annular quantum surface is bounded by $K_0$ with probability $1-\frac\delta 4$. Then for $\phi \sim \LF_{\cC_+, \ell}^{(\alpha, +\infty), (\beta, 0)}$ conditioned on $F_\ell$,  with probability at least $1-\delta - o_\ell(1)$  we have $|(\phi - \frac2\gamma \log \ell, \rho)| < K_0 - \log 2$ and $\cA_{\phi - \frac2\gamma \log \ell}([0,1] \times [0,2\pi]) < K_0$. We are done. 


\end{proof}

\begin{proof}[Proof of Proposition~\ref{prop-sphere-given-F}]
The result is immediate from Lemma~\ref{lem-sphere-given-E} and Proposition~\ref{prop-F|E}.
\end{proof}

\subsection{Proofs of Theorems~\ref{thm-sphere-mot} and~\ref{thm-main}}
\label{subsec-proofs-main}

\begin{proof}[{Proof of Theorem~\ref{thm-sphere-mot}}]
Let $(L_t^\infty, Z_t^\infty)_{[0,\tau^\infty]}$ have the law of the Brownian process described in Theorem~\ref{thm-sphere-mot}.

\noindent \textbf{Step 1: Constructing a pair $(\widetilde \phi^\infty, \widetilde \eta^\infty)$ with boundary length process $(L_t^\infty, Z_t^\infty)$.} 
     For $x >0$ let $\tau_x^\infty$ be the first time $L_t^\infty$ hits $x$ (or, if no such time exists, $\tau_x^\infty = \infty$). 
    For each $\ell$ of the form $2^{-n}$ such that $\tau_\ell^\infty \neq \infty$, by Theorem~\ref{thm-mot-disk-finite} a.s.\ there is a corresponding SLE-decorated quantum surface $\cD_\ell^\infty$ associated to the process $(L_{t - \tau_\ell^\infty}, Z_{t - \tau_\ell^\infty})_{[0, \tau^\infty - \tau^\infty_{\ell}]}$, and the $\cD_\ell^\infty$ are consistent in the sense that for $\ell' < \ell$ the decorated quantum surface $\cD_{\ell}^\infty$ arises as a sub-surface of $\cD_{\ell'}^\infty$. Thus by the Kolmogorov extension theorem there is a curve-decorated quantum surface $(\cC, \widetilde \phi^\infty, \widetilde \eta^\infty, -\infty, +\infty)$ such that for all $\ell = 2^{-n}$ such that {$\tau_\ell^\infty  \neq \infty$}, we have $\cD_\ell^\infty = (\widetilde \eta^\infty([\tau_\ell^\infty, \tau^\infty]), \widetilde\phi^\infty, \widetilde\eta^\infty(\cdot + \tau_\ell^\infty)|_{[0, \tau^\infty - \tau_\ell^\infty]}, \widetilde \eta^\infty(\tau_\ell^\infty), +\infty)$. 

\medskip 

Let $\phi' \sim \cL$ as in Proposition~\ref{prop-sphere-given-F}, so $(\cC, \phi', -\infty, +\infty)/{\sim_\gamma}$ has the law of $\mathcal M_2^\mathrm{sph}(\alpha)$ conditioned to have quantum area greater than 1. Independently let $\eta'$ be whole-plane $\SLE_\kappa$ from $-\infty$ to $+\infty$ in $\cC$. 

\noindent \textbf{Step 2: $(\phi', \eta')$ is the $\ell \to 0$ limit of $\LF_{\cC_+, \ell}^{(\alpha, +\infty), (\beta, 0)}$ decorated by independent radial SLE.} 
    By Proposition~\ref{prop-sphere-given-F},
 for $\phi^\ell \sim \LF_{\cC_+, \ell}^{(\alpha, +\infty), (\beta, 0)}$ conditioned on $F_\ell$, with $\sigma_\ell > 0$ satisfying $\cA_{\phi^\ell}(\cC_+ + \sigma_\ell) = \frac12$ and $\widetilde \phi^\ell = \phi^\ell(\cdot + \sigma_\ell)$, for any $U \subset \mathcal C$ bounded away from {$-\infty$}, as $\ell \to 0$ the field $\widetilde \phi^\ell|_U$ converges in distribution to  $\phi'|_U$. 
    Note that $\sigma_\ell \to \infty$ in probability as $\ell \to 0$ (e.g.\ by taking $U = [-N, \infty) \times [0,2\pi]$ in the above statement). 

    Next, sample a radial $\SLE_\kappa$ curve $\eta^\ell$ in $(\cC_+, 0,+\infty)$ independently of $\phi^\ell$ and parametrize it by quantum area. Let $\widetilde \eta^\ell = \eta^\ell + \sigma_\ell$, and for each neighborhood $U$ of $+\infty$ bounded away from $-\infty$ define the curve $\widetilde \eta_U: [0,\infty)\to \mathcal C$ by $\widetilde \eta_U    := \widetilde \eta(\cdot + \sigma_U)$ where $\sigma_U$ is the first time $\widetilde \eta$ hits $\ol U$. Since whole-plane $\SLE_\kappa$ is the local  limit of radial $\SLE_\kappa$ as the domain tends to the whole plane, the curve $\widetilde \eta_U$ converges in the topology of uniform convergence on compact sets to  $\eta'_U:= \eta'(\cdot + \sigma'_U)$, where $\sigma_U'$ is the time $\eta'$ first hits $\ol U$.
    
	Thus, in the setup of  Theorem~\ref{thm-mot-disk-finite} with boundary length $\ell$ rather than $1$, conditioned on having quantum area at least 1, as $\ell \to 0$ the field and curve $\widetilde \phi^\ell, \widetilde \eta^\ell$ converge in law to $\phi', \eta'$ above. 

\medskip

\noindent \textbf{Step 3: Showing $(\cC,\widetilde \phi^\infty, \widetilde \eta^\infty, -\infty, +\infty)/{\sim_\gamma} \stackrel d= (\cC, \phi', \eta', -\infty, +\infty)/{\sim_\gamma}$.}
For $\ell$ of the form $2^{-n}$, let  $(X_t^\ell, Y_t^\ell)$ be the boundary length process associated to $\widetilde \phi^\ell, \widetilde \eta^\ell$ defined above, and let $(L_t^\ell, Z_t^\ell) = (X_t^\ell - Y_t^\ell, X_t^\ell+Y_t^\ell)$. 
    Since $\tau_\ell^\infty \to 0$ in probability as $\ell \to 0$, we can couple $(L_{t - \tau_\ell^\infty}^\infty, Z_{t - \tau_\ell^\infty}^\infty)$ to agree with $(L_t^\ell, Z_t^\ell)$ with probability $1-o_\ell(1)$. On this event  $(\widetilde\eta^\infty([\tau_\ell^\infty, T^\infty]), \widetilde\phi^\infty, \widetilde \eta^\infty(\cdot + \tau_\ell^\infty)|_{[0, \tau^\infty - \tau_\ell^\infty]}, \widetilde \eta^\infty(\tau_\ell^\infty), +\infty)/{\sim_\gamma} = (\cC_+ - \sigma_\ell, \widetilde\phi^\ell, \widetilde\eta^\ell, 0, +\infty)/{\sim_\gamma}$; let $f_\ell$ be the conformal map sending $\widetilde\eta^\infty([\tau_\ell^\infty, \tau^\infty])$ to $\cC_+ - \sigma_\ell$ such that $f_\ell(\widetilde\eta^\infty(\tau_\ell^\infty)) = -\sigma_\ell$ and $f_\ell(+\infty) =+\infty$. Since for any $N$ the regions
	$\cC \backslash \widetilde\eta^\infty([\tau_\ell^\infty, \tau^\infty])$ and $\cC \backslash (\cC_+ - \sigma_\ell)$  are subsets of $(-\infty, N) \times [0,2\pi]$ in probability as $\ell \to \infty$, 
	standard conformal distortion estimates give that for every neighborhood $U$ of $+\infty$ bounded away from $-\infty$ we have $\sup_U |f'_\ell  - 1|	 \to 0$ in probability. 
 This implies that there is a coupling of $(\phi', \eta')$ with $(\widetilde\phi^\infty, \widetilde\eta^\infty)$ and a random rotation  $f_\infty: \cC \to \cC$ of the cylinder (i.e.\ conformal map fixing $\pm \infty$ with $\mathrm{Re}\, f_\infty(z) = \mathrm{Re}\, z$ for all $z$) such that $\widetilde\phi^\infty = f_\infty \bullet_\gamma(\phi')$ and $\widetilde\eta^\infty = f_\infty (\eta')$ a.s., completing the step. 

 \medskip \noindent \textbf{Conclusion.} $(\cC, \widetilde \phi^\infty, \widetilde \eta^\infty, -\infty, +\infty)/{\sim_\gamma}$ has the law of the curve-decorated quantum surface of Theorem~\ref{thm-sphere-mot} (Step 3), and its boundary length process is as desired (Step 1).  The measurability claim~\eqref{eq-thm-sph} is immediate from that of Theorem~\ref{thm-mot-disk-finite} and the construction of Step 1. 
\end{proof}

\begin{figure}[ht]
	\centering
	\includegraphics[scale=0.5]{figures/reversibility.pdf}
	\caption{Proof of Theorem~\ref{thm-main}. \textbf{Top:} $(\bbC, \phi, 0, \infty)$ is an embedding of a sample from $\mathcal M_2^\mathrm{sph}(Q-\frac\gamma4)$ and $\eta$ is an independent whole-plane SLE from 0 to $\infty$.  Let $(L, Z)$ be its boundary length process. The grey Brownian motion segments are $(L,Z)_{[0,\eps]}$ and $(L,Z)_{[\tau-\eps, \tau]}$. The map $F$ identifies the red, green and blue Brownian motion segments with the corresponding quantum cells.  \textbf{Right:} $(\wt \phi, \wt \eta)$ is obtained from $(\phi, \eta)$ by inverting the plane and orienting $\wt \eta$ so it is a curve from $0$ to $\infty$. \textbf{Left:} $(\wt L, \wt Z)$ is the time-reversal of $(L, Z)$ (translated to start at $0$). By reversibility of Brownian motion we have $(L, Z) \stackrel d= (\wt L, \wt Z)$.  \textbf{Bottom:} By the reversibility of $F$, the map $F$ identifies the red, green and blue Brownian motion segments with the corresponding quantum cells. Sending $\eps \to 0$, we see $F^\infty ((\wt L, \wt Z)) = (\bbC, \wt \phi, \wt \eta,  \infty,0)/{\sim_\gamma}$, so $(\bbC, \phi, \eta, \infty,0)/{\sim_\gamma} = F^\infty((L,Z)) \stackrel d= F^\infty((\wt L, \wt Z)) = (\bbC, \wt \phi, \wt \eta, \infty,0)/{\sim_\gamma}$. This implies reversibility of whole-plane SLE.}
	\label{fig:reversibility}
\end{figure}

\begin{proof}[{Proof of Theorem~\ref{thm-main}}]
    In the setting of Theorem~\ref{thm-sphere-mot},
      {the decorated quantum surface $(\bbC, \phi, \eta, \infty, 0)/{\sim_\gamma}$ is measurable with respect to $(L, Z)$. Indeed, let $(X_t, Y_t) = (\frac12(L_t + Z_t), \frac12 (L_t - Z_t))$ and let $(t_n)_{n \in \bbZ}$ be an increasing collection of rational times in $(0,\tau)$ such that $\lim_{n \to -\infty} t_n = 0$, $\lim_{n \to \infty} t_n = \tau$, and for each $n$ we have $(X_{t_n} - \inf_{[t_n, t_{n+1}]} X_\cdot)+ (Y_{t_n} - \inf_{[t_n, t_{n+1}]} Y_\cdot) < L_{t_n}$ (i.e., $\eta([t_n,t_{n+1}])$ is simply connected). Then $(\bbC, \phi, \eta, \infty, 0)/{\sim_\gamma}$ is the conformal welding of $\cC_n := F((X_{\cdot + t_n} - X_{t_n}, Y_{\cdot + t_n} - Y_{t_n})|_{[0,t_{n+1}-t_n]})$ for $n \in \bbZ$, where, similarly as in Figure~\ref{fig-lengths} (right),  the first marked point of $\cC_{n+1}$ is identified with the second marked point of $\cC_n$, and the two boundary arcs of $\cC_{n+1}$ adjacent to its first marked point are conformally welded according to quantum length to the boundary of the conformal welding of $\bigcup_{i \leq n} \cC_i$.} 
    Let $F^\infty$ be the map sending the process $(L_t, Z_t)_{[0,\tau]}$ to the decorated quantum surface $(\bbC, \phi, \eta, \infty, 0)/{\sim_\gamma}$.  

See Figure~\ref{fig:reversibility}.
    Let $(\bbC, \phi, \infty, 0)$ be an embedding of a sample from $\mathcal M^\mathrm{sph}_2(Q - \frac\gamma4)$ and let $\eta$ be an independent whole-plane SLE. Let $(L_t, Z_t)$ be the boundary length process associated to $(\bbC, \phi, \eta, \infty, 0)$ as in Theorem~\ref{thm-sphere-mot},  {and let $\tau$ be its random duration}. Let $\mathrm{Inv}(z) = z^{-1}$, let $\widetilde \phi = \mathrm{Inv} \bullet_\gamma \phi$, let $\widetilde \eta$ be the time-reversal of $\mathrm{Inv} \circ \eta$ {(so $\widetilde \eta$ is also a curve from $0$ to $\infty$)}, and let $(\widetilde L_t, \widetilde Z_t)  { := (L_{\tau - t}, Z_{\tau - t} - Z_\tau)}$ be the time-reversal of $(L_t, Z_t)$. Let $\mathfrak S = (\bbC, \phi, \eta, \infty, 0)$ and $\widetilde {\mathfrak S} = (\bbC, \widetilde \phi, \widetilde \eta, \infty, 0)$.  {In the next paragraph we will show that $\mathfrak S \stackrel d= \wt {\mathfrak S}$; this is the crux of the argument.}
    
     {By definition $F^\infty((L, Z)) = \mathfrak S$. 
     We pick $\e,n>0$ and equally divide $(\e,\tau-\e)$ into $n$ intervals $I_1,...,I_n$. We restrict to the event $E_n$ that $\eta({I_k})$ is simply connected for $k=1,...,n$. Let $\cC_k = (\eta(I_k),\phi,\eta|_{I_k})/{\sim_\gamma}$ and $\wt \cC_k =  (\wt\eta(I_k),\wt\phi,\wt\eta|_{I_k})/{\sim_\gamma}$, so by definition of $(\wt \phi, \wt \eta)$ the decorated quantum surface $\wt\cC_{n+1-k}$ agrees with $\cC_k$ with its curve reversed. 
     	For an interval $I = [a,b]$, define $F'((L, Z)|_I) := F((X_{\cdot + a} - X_a, Y_{\cdot + a} - Y_a)|_{[0, b-a]})$ where $(X_\cdot, Y_\cdot) = (\frac12 (L_\cdot+Z_\cdot),\frac12 (L_\cdot-Z_\cdot))$. 
     	By  Theorem~\ref{thm-sphere-mot} $F'((L,Z)|_{I_k}) = \cC_k$, so by reversibility of $F$ (Lemma~\ref{lem-F-reversible}) and the fact that $(L, Z)|_{I_k}$ and $(\wt L, \wt Z - Z_\tau)|_{I_{n+1-k}}$ differ by time-reversal, we have $F((\wt L,\wt Z)|_{I_{n+1-k}}) = \wt\cC_{n+1-k}$. Consequently, on the event $E_n$, the boundary length process of $(\wt \phi, \wt \eta)$ restricted to $(\eps, \tau - \eps)$ agrees with $(\frac12 (\wt L + \wt Z), \frac12 (\wt L - \wt Z))|_{[\eps, \tau-\eps]}$ up to additive constant.
     	Therefore, by first sending $n\to\infty$ and then $\e\to0$, we see that $F^\infty((\widetilde L, \widetilde Z)) = \widetilde {\mathfrak S}$ a.s.. } 
    The reversibility of Brownian motion yields $(L_t, Z_t) \stackrel d= (\widetilde L_t, \widetilde Z_t)$, and combining with $(F^\infty((L, Z)), F^\infty ((\wt L, \wt Z)) = (\mathfrak S, \wt {\mathfrak S})$, we conclude $\mathfrak S \stackrel d= \widetilde {\mathfrak S}$. 
    
    Let $r>0$ be such that $\cA_\phi(r \bbD) = \cA_\phi( \bbC \backslash r \bbD)$, let $\theta$ be uniformly sampled from $[0,2\pi)$ independently of $(\phi, \eta)$, define $f(z) = r^{-1} e^{i\theta} z$, and set $\phi_0 = f \bullet_\gamma \phi$ and $\eta_0 = f \circ \eta$. Likewise define $\widetilde \phi_0, \widetilde \eta_0$ by applying the same embedding procedure for $\widetilde \phi, \widetilde \eta$. Since $\mathfrak S \stackrel d= \widetilde {\mathfrak S}$ we have $(\phi_0, \eta_0) \stackrel d= (\widetilde \phi_0, \widetilde \eta_0)$. Since $\phi$ and $\eta$ are independent, and whole-plane SLE is invariant in law under dilations and rotations of the plane, the law of $\eta_0$ is whole-plane SLE. Likewise  $\widetilde \eta_0$ has the law of the time-reversal of whole-plane SLE after applying $\mathrm{Inv}$. The statement $\eta_0 \stackrel d= \widetilde \eta_0$ is thus the desired reversibility of whole-plane SLE for $\kappa > 8$.
\end{proof}

\section{Open problems}\label{sec-open}

{
The convergence of lattice statistical physics models to SLE was a primary reason to expect the reversibility of chordal $\SLE_\kappa$ for  $\kappa \leq 8$. Conversely, Theorem~\ref{thm-main} suggests the following question. 

\begin{prob}
    Find a  lattice statistical physics model whose scaling limit is whole-plane $\SLE_\kappa$ for some $\kappa > 8$. 
\end{prob}

Questions of this sort are sometimes easier when the underlying lattice is random, i.e., is a random planar map. Some random planar maps decorated by statistical physics models can be encoded by a pair of trees, which in turn may be described by a random walk on the 2D lattice. If this random walk converges in the scaling limit to Brownian motion with covariance given by~\eqref{eqn-thm:mot-chordal}, then we say the corresponding decorated random planar map converges in the \emph{peanosphere topology} to $(\gamma = \frac4{\sqrt\kappa})$-LQG decorated by space-filling $\SLE_{\kappa}$.
In the case when the $\SLE_{\kappa}$ is a space-filling loop in $\hat \bbC$ from $\infty$ to $\infty$, such convergences are known for random planar maps decorated by bipolar orientations ($\kappa = 12$) \cite{KMSW19}, Schnyder woods ($\kappa = 16$) \cite{LSW17}, or a variant of spanning trees ($\kappa > 8$) \cite{GKMW18}; see the survey \cite{ghs-mating-survey} for examples where $\kappa \leq 8$. The next problem asks for such a result for whole-plane $\SLE_\kappa$ where $\kappa > 8$. 
}
\begin{prob}\label{prob-model}
   Exhibit a random planar map decorated by a statistical physics model which can be encoded by a random walk converging in the limit to $(X_t, Y_t)$ defined in Theorem~\ref{thm-mot-disk-finite} or in Theorem~\ref{thm-sphere-mot}. In other words, find a random planar map model which converges in the peanosphere sense to LQG decorated by radial or whole-plane SLE. 
\end{prob}

 {One of the most natural variants of  $\SLE_\kappa$ is the $\SLE_\kappa(\underline\rho)$ process \cite{LSW03, Dub09, MS16a}; other important variants include multiple SLE~\cite{bauer2005multiple,kozdron2006configurational,dubedat2007commutation} and the conformal loop ensemble~\cite{SheffieldCLE,Sheffield-Werner-CLE}.} For $\kappa\in (0,8]$, the time reversal of chordal $\SLE_\kappa(\underline\rho)$ has been solved when the sum of the weights is larger than $(-2)\vee(\frac{\kappa}{2}-4)$~\cite{MS16b,Zhan21,Yu22}, while~\cite[Theorem 1.18]{MS17} gives a criterion for the reversibility of $\SLE_\kappa(\rho^-;\rho^+)$ curves when $\kappa>8$. 
On the whole plane side, the most natural variant of whole-plane $\SLE_\kappa$ is  whole-plane $\SLE_\kappa(\rho)$ for $\rho > -2$, which agrees  with whole-plane $\SLE_\kappa$ when $\rho=0$ (see e.g.~\cite[Section 2.1.3]{MS17}). Miller and Sheffield showed that when $\kappa \in (0,4]$ and $\rho > -2$, or $\kappa \in (4, 8]$ and $\rho \geq \frac\kappa2-4$,  whole-plane $\SLE_\kappa(\rho)$ is reversible \cite[Theorem 1.20]{MS17}. They also show that when $\kappa > 8$ and $\rho \geq \frac\kappa2 - 4$  whole-plane $\SLE_\kappa(\rho)$ is not reversible \cite[Remark 1.21]{MS17}. They do not treat the regime where $\kappa > 4$ and $\rho \in (-2, \frac\kappa2-4)$ because it is not as natural in the imaginary geometry framework, see \cite[Remark 1.22]{MS17}. On the other hand, Theorem~\ref{thm-main} gives reversibility when  $\kappa > 8$ and $\rho = 0$ even though it falls into this regime, so there is still hope for reversibility in this range.  

\begin{prob}\label{prob-bigger-8}
    When $\kappa > 4$, for which $\rho \in (-2, \frac\kappa2-4)$ is $\SLE_\kappa(\rho)$ reversible? 
\end{prob}

A further generalization of whole-plane $\SLE_\kappa(\rho)$ can be obtained by adding a constant drift term to the driving function. 
Zhan~\cite{zhan-whole-plane} showed that when $\kappa \in (0,4]$, $\rho = 0$ and any constant drift is chosen, the curve is reversible. Miller and Sheffield~\cite[Theorem 1.20]{MS17} showed that when $\kappa \in (0,4]$ and $\rho > -2$, or $\kappa \in (4,8]$ and $\rho \geq \frac\kappa2-4$, for any chosen drift the curve is reversible. 

\begin{prob}\label{prob-drift}
When $\kappa > 4$ and $\rho \in (-2, \frac\kappa2-4)$, what choices of drift coefficient give a reversible curve? 
\end{prob}

The statement of Theorem~\ref{thm-main}  involves only SLE, but our arguments depend on couplings with LQG. 

\begin{prob}\label{prob-not-mot}
    Find a proof of Theorem~\ref{thm-main} not using mating-of-trees.
\end{prob}
It seems likely that a solution to Problem~\ref{prob-not-mot} would represent a significant step towards solving Problems~\ref{prob-bigger-8} and~\ref{prob-drift}.

\bibliographystyle{alpha}
\bibliography{theta}

\begin{thebibliography}{GKMW18}

\bibitem[AG21]{AG21}
Morris Ang and Ewain Gwynne.
\newblock {Liouville quantum gravity surfaces with boundary as matings of
  trees}.
\newblock {\em Annales de l'Institut Henri Poincar\'e, Probabilit\'es et
  Statistiques}, 57(1):1 -- 53, 2021.

\bibitem[AHS23]{AHS21}
Morris Ang, Nina Holden, and Xin Sun.
\newblock {Integrability of SLE via conformal welding of random surfaces}.
\newblock {\em Communications on Pure and Applied Mathematics}, pages 1--57,
  2023.

\bibitem[Ang23]{Ang23}
Morris Ang.
\newblock Liouville conformal field theory and the quantum zipper.
\newblock {\em arXiv preprint arXiv:2301.13200}, 2023.

\bibitem[ARS22]{ars-annulus}
Morris {Ang}, Guillaume {Remy}, and Xin {Sun}.
\newblock {The moduli of annuli in random conformal geometry}.
\newblock {\em ArXiv e-prints}, March 2022.

\bibitem[ARS23]{ARS21}
Morris Ang, Guillaume Remy, and Xin Sun.
\newblock {FZZ formula of boundary Liouville CFT via conformal welding}.
\newblock {\em Journal of the European Mathematical Society}, pages 1--58,
  2023.

\bibitem[ARSZ23]{ARSZ23}
Morris Ang, Guillaume Remy, Xin Sun, and Tunan Zhu.
\newblock Derivation of all structure constants for boundary {L}iouville {CFT}.
\newblock {\em arXiv preprint arXiv:2305.18266}, 2023.

\bibitem[AS21]{AS21}
Morris Ang and Xin Sun.
\newblock Integrability of the conformal loop ensemble.
\newblock {\em arXiv preprint arXiv:2107.01788}, 2021.

\bibitem[ASY22]{ASY22}
Morris Ang, Xin Sun, and Pu~Yu.
\newblock Quantum triangles and imaginary geometry flow lines.
\newblock {\em arXiv preprint arXiv:}, 2022.

\bibitem[BBK05]{bauer2005multiple}
Michel Bauer, Denis Bernard, and Kalle Kyt{\"o}l{\"a}.
\newblock Multiple {S}chramm--{L}oewner evolutions and statistical mechanics
  martingales.
\newblock {\em Journal of statistical physics}, 120:1125--1163, 2005.

\bibitem[BN]{bn-sle-notes}
{N}.\ {B}erestycki and {J}.{R}.\ {N}orris.
\newblock {L}ectures on {S}chramm-{L}oewner {E}volution.
\newblock Available at \url{http://www.statslab.cam.ac.uk/~james/Lectures/}.

\bibitem[DKRV16]{DKRV16}
Fran{\c{c}}ois David, Antti Kupiainen, R{\'e}mi Rhodes, and Vincent Vargas.
\newblock Liouville quantum gravity on the {R}iemann sphere.
\newblock {\em Communications in Mathematical Physics}, 342(3):869--907, 2016.

\bibitem[DMS21]{DMS14}
Bertrand Duplantier, Jason Miller, and Scott Sheffield.
\newblock Liouville quantum gravity as a mating of trees.
\newblock {\em Ast\'{e}risque}, 427, 2021.

\bibitem[DS11]{DS11}
Bertrand Duplantier and Scott Sheffield.
\newblock Liouville quantum gravity and {KPZ}.
\newblock {\em Inventiones mathematicae}, 185(2):333--393, 2011.

\bibitem[Dub07]{dubedat2007commutation}
Julien Dub{\'e}dat.
\newblock Commutation relations for {S}chramm-{L}oewner evolutions.
\newblock {\em Communications on Pure and Applied Mathematics: A Journal Issued
  by the Courant Institute of Mathematical Sciences}, 60(12):1792--1847, 2007.

\bibitem[Dub09]{Dub09}
J.~Dub\'{e}dat.
\newblock Duality of {S}chramm-{L}oewner {E}volutions.
\newblock {\em Ann. Sci. \'{E}c. Norm. Sup\'{e}r}, 42(5), 2009.

\bibitem[Dur63]{annular-distortion}
Peter~L Duren.
\newblock Distortion in certain conformal mappings of an annulus.
\newblock {\em Michigan Mathematical Journal}, 10(4):431--441, 1963.

\bibitem[GHMS17]{GHMS15}
Ewain Gwynne, Nina Holden, Jason Miller, and Xin Sun.
\newblock Brownian motion correlation in the peanosphere for $\kappa'> 8$.
\newblock In {\em Annales de l'Institut Henri Poincar{\'e}, Probabilit{\'e}s et
  Statistiques}, volume~53, pages 1866--1889. Institut Henri Poincar{\'e},
  2017.

\bibitem[GHS23]{ghs-mating-survey}
Ewain Gwynne, Nina Holden, and Xin Sun.
\newblock Mating of trees for random planar maps and {L}iouville quantum
  gravity: a survey.
\newblock In {\em Topics in statistical mechanics}, volume~59 of {\em Panor.
  Synth\`eses}, pages 41--120. Soc. Math. France, Paris, 2023.

\bibitem[GKMW18]{GKMW18}
Ewain Gwynne, Adrien Kassel, Jason Miller, and David~B Wilson.
\newblock Active spanning trees with bending energy on planar maps and
  {SLE}-decorated {L}iouville quantum gravity for $\kappa> 8$.
\newblock {\em Communications in Mathematical Physics}, 358:1065--1115, 2018.

\bibitem[GMS18]{GMS18}
Ewain Gwynne, Jason Miller, and Xin Sun.
\newblock {Almost sure multifractal spectrum of {S}chramm-{L}oewner evolution}.
\newblock {\em Duke Mathematical Journal}, 167(6):1099 -- 1237, 2018.

\bibitem[Hos01]{hosomichi}
Kazuo Hosomichi.
\newblock Bulk-boundary propagator in {L}iouville theory on a disc.
\newblock {\em Journal of High Energy Physics}, 2001(11):044, 2001.

\bibitem[HRV18]{HRV-disk}
Yichao Huang, R\'{e}mi Rhodes, and Vincent Vargas.
\newblock Liouville quantum gravity on the unit disk.
\newblock {\em Ann. Inst. Henri Poincar\'{e} Probab. Stat.}, 54(3):1694--1730,
  2018.

\bibitem[KL06]{kozdron2006configurational}
Michael~J Kozdron and Gregory~F Lawler.
\newblock The configurational measure on mutually avoiding {SLE} paths.
\newblock {\em arXiv preprint math/0605159}, 2006.

\bibitem[KMS23]{KMS-removability}
Konstantinos Kavvadias, Jason Miller, and Lukas Schoug.
\newblock {Conformal removability of non-simple Schramm-Loewner evolutions}.
\newblock {\em arXiv preprint arXiv:2302.10857}, 2023.

\bibitem[KMSW19]{KMSW19}
Richard Kenyon, Jason Miller, Scott Sheffield, and David~B Wilson.
\newblock Bipolar orientations on planar maps and {SLE}$_{12}$.
\newblock {\em The Annals of Probability}, 47(3):1240--1269, 2019.

\bibitem[Law08]{Law08}
Gregory~F Lawler.
\newblock {\em Conformally invariant processes in the plane}.
\newblock Number 114. American Mathematical Soc., 2008.

\bibitem[Law13]{Lawler11continuity}
Gregory~F. Lawler.
\newblock Continuity of radial and two-sided radial sle at the terminal point.
\newblock In {\em In the tradition of {A}hlfors-{B}ers. {VI}}, volume 590 of
  {\em Contemp. Math.}, pages 101--124. Amer. Math. Soc., Providence, RI, 2013.

\bibitem[LSW03]{LSW03}
Gregory Lawler, Oded Schramm, and Wendelin Werner.
\newblock Conformal restriction: the chordal case.
\newblock {\em Journal of the American Mathematical Society}, 16(4):917--955,
  2003.

\bibitem[LSW04]{lsw-lerw-ust}
Gregory~F. Lawler, Oded Schramm, and Wendelin Werner.
\newblock Conformal invariance of planar loop-erased random walks and uniform
  spanning trees.
\newblock {\em Ann. Probab.}, 32(1B):939--995, 2004.

\bibitem[LSW24]{LSW17}
Yiting Li, Xin Sun, and Samuel Watson.
\newblock {Schnyder woods, SLE$_{16}$, and Liouville quantum gravity}.
\newblock {\em Transactions of the American Mathematical Society}, 2024.

\bibitem[MS16a]{MS16a}
J.~Miller and S.~Sheffield.
\newblock Imaginary {G}eometry {I}: Interacting {SLE}s.
\newblock {\em Probability Theory and Related Fields}, 164(3-4):553--705, 2016.

\bibitem[MS16b]{MS16b}
Jason Miller and Scott Sheffield.
\newblock {Imaginary geometry II: Reversibility of
  {SLE}$_\kappa(\rho_1;\rho_2)$ for $\kappa\in (0, 4) $}.
\newblock {\em The Annals of Probability}, 44(3):1647--1722, 2016.

\bibitem[MS16c]{ig3}
Jason Miller and Scott Sheffield.
\newblock Imaginary geometry {III}: reversibility of {$\mathrm{SLE}_\kappa$}
  for {$\kappa\in(4,8)$}.
\newblock {\em Ann. of Math. (2)}, 184(2):455--486, 2016.

\bibitem[MS16d]{ms-qle}
Jason Miller and Scott Sheffield.
\newblock {Quantum Loewner evolution}.
\newblock {\em Duke Mathematical Journal}, 165(17):3241 -- 3378, 2016.

\bibitem[MS17]{MS17}
Jason Miller and Scott Sheffield.
\newblock {Imaginary geometry IV: interior rays, whole-plane reversibility, and
  space-filling trees}.
\newblock {\em Probability Theory and Related Fields}, 169(3):729--869, 2017.

\bibitem[MS19]{MS19}
Jason Miller and Scott Sheffield.
\newblock Liouville quantum gravity spheres as matings of finite-diameter
  trees.
\newblock {\em Annales de l'Institut Henri Poincar{\'e}, Probabilit{\'e}s et
  Statistiques}, 55(3):1712--1750, 2019.

\bibitem[RS05]{RS05}
S.~Rohde and O.~Schramm.
\newblock Basic properties of {SLE}.
\newblock {\em Ann. of Math.}, 161(2), 2005.

\bibitem[RV10]{RV10}
Raoul Robert and Vincent Vargas.
\newblock Gaussian multiplicative chaos revisited.
\newblock {\em The Annals of Probability}, 38(2):605--631, 2010.

\bibitem[RZ16]{zhanrohde2016}
Steffen Rohde and Dapeng Zhan.
\newblock Backward {SLE} and the symmetry of the welding.
\newblock {\em Probability Theory and Related Fields}, 164(3-4):815--863, 2016.

\bibitem[RZ20]{RZ20b}
Guillaume Remy and Tunan Zhu.
\newblock {Integrability of boundary Liouville conformal field theory}.
\newblock {\em arXiv preprint arXiv:2002.05625}, 2020.

\bibitem[Sch00]{Sch00}
Oded Schramm.
\newblock Scaling limits of loop-erased random walks and uniform spanning
  trees.
\newblock {\em Israel Journal of Mathematics}, 118(1):221--288, 2000.

\bibitem[She09]{SheffieldCLE}
Scott Sheffield.
\newblock Exploration trees and conformal loop ensembles.
\newblock {\em Duke Math. J.}, 147(1):79--129, 2009.

\bibitem[She16]{She16a}
Scott Sheffield.
\newblock Conformal weldings of random surfaces: {SLE} and the quantum gravity
  zipper.
\newblock {\em The Annals of Probability}, 44(5):3474--3545, 2016.

\bibitem[Smi01]{smirnov-cardy}
Stanislav Smirnov.
\newblock Critical percolation in the plane: conformal invariance, {C}ardy's
  formula, scaling limits.
\newblock {\em C. R. Acad. Sci. Paris S\'er. I Math.}, 333(3):239--244, 2001.

\bibitem[Smi10]{smirnov-ising}
Stanislav Smirnov.
\newblock Conformal invariance in random cluster models. {I}. {H}olomorphic
  fermions in the {I}sing model.
\newblock {\em Ann. of Math. (2)}, 172(2):1435--1467, 2010.

\bibitem[SS09]{ss-dgff}
Oded Schramm and Scott Sheffield.
\newblock Contour lines of the two-dimensional discrete {G}aussian free field.
\newblock {\em Acta Math.}, 202(1):21--137, 2009.

\bibitem[SW12]{Sheffield-Werner-CLE}
Scott Sheffield and Wendelin Werner.
\newblock Conformal loop ensembles: the {M}arkovian characterization and the
  loop-soup construction.
\newblock {\em Ann. of Math. (2)}, 176(3):1827--1917, 2012.

\bibitem[SW16]{SW16}
Scott Sheffield and Menglu Wang.
\newblock {Field-measure correspondence in Liouville quantum gravity almost
  surely commutes with all conformal maps simultaneously}.
\newblock {\em arXiv preprint arXiv:1605.06171}, 2016.

\bibitem[VW20]{vw-loewner-kufarev}
Fredrik Viklund and Yilin Wang.
\newblock The {L}oewner-{K}ufarev {E}nergy and {F}oliations by
  {W}eil-{P}etersson {Q}uasicircles.
\newblock {\em ArXiv e-prints}, 2020.

\bibitem[Yu23]{Yu22}
Pu~Yu.
\newblock {Time-reversal of multiple-force-point chordal
  $\mathrm{SLE}_\kappa(\underline{\rho})$}.
\newblock {\em Electronic Journal of Probability}, to appear, 2023.

\bibitem[Zha08a]{zhan2008duality}
Dapeng Zhan.
\newblock Duality of chordal {SLE}.
\newblock {\em Inventiones mathematicae}, 174(2):309--353, 2008.

\bibitem[Zha08b]{zhan-chordal}
Dapeng Zhan.
\newblock Reversibility of chordal {SLE}.
\newblock {\em Ann. Probab.}, 36(4):1472--1494, 2008.

\bibitem[Zha15]{zhan-whole-plane}
Dapeng Zhan.
\newblock Reversibility of whole-plane {SLE}.
\newblock {\em Probability Theory and Related Fields}, 161(3-4):561--618, 2015.

\bibitem[Zha21]{Zhan21}
Dapeng Zhan.
\newblock Boundary {G}reen's functions and {M}inkowski content measure of
  multi-force-point {SLE}$_\kappa(\underline\rho) $.
\newblock {\em arXiv preprint arXiv:2106.12670}, 2021.

\end{thebibliography}

\end{document}